\documentclass[twocolumn]{svjour3}          

\smartqed

\usepackage[english]{babel} 
\usepackage[T1]{fontenc} 
\usepackage[ansinew]{inputenc} 
\usepackage{lmodern} 
\usepackage{hyperref}
\usepackage{graphicx}
\usepackage{listings}
\usepackage{subfig}
\usepackage{amsmath}
\usepackage{amsfonts}
\usepackage{amssymb}
\usepackage{booktabs} 
\usepackage{etoolbox}
\usepackage{multirow}
\usepackage{float}
\usepackage{natbib}

\journalname{}

\begin{document}

\title{Revisiting topology optimization with buckling constraints}

\titlerunning{Revisiting topology optimization with buckling constraints}

\author{Federico Ferrari \and
        Ole Sigmund}

\authorrunning{F. Ferrari \and O. Sigmund}

\institute{Department of Mechanical Engineering \\
           Technical University of Denmark \\
           Nils Koppels All\'{e} 404, 2800 Kongens Lyngby, Denmark \\
           \email{feferr@mek.dtu.dk, sigmund@mek.dtu.dk} \\
           }

\date{Received: date / Accepted: date}

\maketitle

\begin{abstract}
 We review some features of topology optimization with a lower bound on the critical load factor, as computed by linearized buckling analysis. The change of the optimized design, the competition between stiffness and stability requirements and the activation of several buckling modes, depending on the value of such lower bound, are studied. We also discuss some specific issues which are of particular interest for this problem, as the use of non-conforming finite elements for the analysis, the use of inconsistent sensitivities in the optimization and the replacement of the single eigenvalue constraints with an aggregated measure. We discuss the influence of these practices on the optimization result, giving some recommendations.
\end{abstract}

\section{Introduction}
 \label{Sec:Intro}

Stability and buckling have attracted considerable attention since early times of structural optimization, due to their importance in the design of structural elements. Moreover, the optimal design according to weight or compliance minimization may naturally lead to structural configurations showing poor stability \citep{rozvany_96a,book:makoto-ohsaki2007}.

Optimal design with respect to buckling has been thoroughly studied for beam models, where even a closed form expression for the sufficient optimality condition can be found \citep{olhoff-rasmussen_77a,cox-overton_92a,kirmser-hu_95a,cox-mccarthy_98a}. Such a goal has turned out to be much more challenging already for plate models \citep{armand-lodier_78a}. In this case better results have been obtained by formulating a reinforcement problem \citep{simitses_73a,haftka-prasad_81a}, rather than allowing for only a continuously varying thickness \citep{frauenthal_72a}.

Much research effort has been devoted to the design of trusses, where stability was at first imposed on the local level, constraining the maximum stress or displacement of individual members \citep{berke_70a,achtziger_99a}. Then, several methods have been proposed to account for overall stability of structures \citep{khot-etal_76a,szyszkowsky-watson_88a,rozvany_96a}. The interaction between local and global buckling modes, of particular concern for building frames stiffened with bracing systems, was studied by \citet{hall-etal_88a}.

Relatively few works have appeared concerning buckling in the optimal design of continuum models. In this setting the problem becomes much more complicated, both for the less intuitive definition of the buckling mechanism and because of several issues hampering the optimization process. 

On one hand, it is difficult to identify the optimized design with a grid of tension and compression members and therefore the buckling mechanism becomes less intuitive. On the other hand, typical issues encountered in eigenvalue optimization, such as eigenvalue multiplicity \citep{seyranian_94a}, artificial modes \citep{neves-etal_95a,neves-etal_02a} and existence of many local minima are encountered, posing serious convergence issues. A large number of buckling modes are likely to become active and therefore must be considered \citep{bruyneel-etal_08a,dunning-etal_16a}, making the problem also very challenging from a computational point of view.

Homogenization--based topology optimization aimed at the maximization of the linearized buckling load was first addressed by \citet{neves-etal_95a}, while \citet{rodrigues-etal_95a} provided a comprehensive study of the optimality conditions for both single and multiple eigenvalues. This latter situation was then accounted for by \citet{folgado-rodrigues_98a}, while solving the reinforcement problem for a plate.

Topology Optimization (TO) is constantly spreading to more fields of structural engineering as a preliminary design tool, also for large--scale, complex structures \citep{aage-etal_17a}. Due to the specific and often counter--acting character of buckling design, as opposed to compliance or stress designs, it is apparent that the former must be carefully taken into account in the optimization of such complex structures.

Therefore a renewed interest in topology design with regards to buckling is raising. Some works have focused on alleviating the issues due to the clashing character of stiffness and stability \citep{gao-ma_15a,gao-etal_17a}, or on the use of effective iterative methods for solving the large eigenvalue problem \citep{dunning-etal_16a,bian-feng_17a}. Recently, \citet{thomsen-etal_18a} introduced a method for the design of periodic microstructures with respect to multi--scale buckling conditions, laying the foundations for plenty of future applications within multi--scale structural and material design. Other recent works and applications can be found where we wish to mention \citet{zhou_04a,lund_09a,bochenek-tajs_15a,cheng-xu_16,chin-kennedy_16a}.

The goal of this note is to provide a discussion about the influence of stability requirements, here posed in the shape of a constraint on the linearized fundamental buckling load factor, in a minimum compliance TO problem. A simple but illustrative example is studied, focusing on the modification of the design and its performance as this constraint becomes more demanding, and on the progressive activation of more buckling modes with either global and local character.

We also discuss the use of aggregation functions for replacing the constraint on the, possibly non--smooth, lowest eigenvalue with a smooth approximation \citep{chen-etal_04a}. Aggregation functions have been extensively used in the context of TO, mainly for stress constrained problems \citep{yang-chen_96a,duysinx-sigmund_98a,le-etal_10a,lee-etal_12a,lee-etal_16a,verbart-etal_17a}. Relatively few works discuss aggregation of eigenvalues, for example in the context of dynamic problems \citep{manh-etal_11a,torii-faria_17a}, and we have been able to locate only the work from \citet{chin-kennedy_16a} concerning buckling constraints.

The influence of some other issues, such as the use of non--conforming finite elements for the buckling analysis, as opposed to the popular conforming four node elements, or the adoption of a simplified but inconsistent sensitivity expression, is also investigated. An ultimate and exhaustive discussion of these topics is clearly outside the aims of this paper. Nonetheless, we believe that the insight provided here might serve as useful guidelines to researchers working in the field.

Finally, it is fair to point out that the linearized analysis addressed here has always raised considerable criticism in the engineering community as, for many structures, it gives an overly simplified description of a complex phenomenon as buckling \citep{kerr-soifer_68a,brantman_77a}. Therefore, in all but very few cases, the estimation of buckling should be carried out by the more appropriate, yet much more computationally expensive, procedure based on non--linear equilibrium \citep{bathe-dvorkin_83a}. Non--linear buckling has been included in \cite{wu-arora_1988} and TO based on it has been carried out by several other authors, \citep[c.f.][]{rahmatalla-swan_03a,kemmler-etal_05a,lindgaard-dahl_13a}. A critical review and comparison of the different approaches to buckling estimation for use in optimal design is given in the recent work of \citet{pedersen-pedersen_18b}.

However, the linearized pre--buckling analysis is still very popular in TO, mainly due to simpler implementation and its computational cheapness. Therefore, we consider a discussion on it to be worthwhile.

The outline of the paper is as follows: in \autoref{Sec:Setting} we introduce the basic settings and recall the theoretical elements used in the further discussion. In \autoref{Sec:NumExp} we present a numerical application, discussing issues touched upon earlier. Finally, conclusions are drawn and some recommendations are given in \autoref{Sec:ConcludingDiscussion}.

\section{Setting and optimization problem}
\label{Sec:Setting}

Let us consider a discretized mechanical system, re\-pre\-sented by $m$ finite elements and $n$ Degrees Of Freedom (DOFs). We address linearized (eigenvalue) buckling analysis \citep{book:cook-etal01}. This requires the selection of a reference load $\mathbf{f}_{0} \in \mathbb{R}^{n}$ and the solution of the linear system
\begin{equation}
 \label{eq:LinearizedEquilibriumEquation}
  \mathbf{K}\mathbf{u}_{0} = \mathbf{f}_{0}
\end{equation}
where $\mathbf{K}\in\mathbb{R}^{n\times n}$ is the linear, symmetric and positive definite stiffness matrix and $\mathbf{u}_{0}\in\mathbb{R}^{n}$ is the equilibrium displacement vector. Given this, the symmetric but indefinite stress stiffness matrix $\mathbf{K}_{\sigma}\left( \mathbf{u}_{0} \right)\in\mathbb{R}^{n\times n}$ can be set up and the generalized eigenvalue problem to be solved is
\begin{equation}
 \label{eq:GEPforMu}
  \left[ \mathbf{K} + \lambda\mathbf{K}_{\sigma}\left( \mathbf{u}_{0} \right) \right] \boldsymbol{\varphi} = \mathbf{0} \: ,
  \qquad \boldsymbol{\varphi} \neq \mathbf{0}
\end{equation}

The eigenpairs $\left( \lambda_{i}, \boldsymbol{\varphi}_{i} \right)$, $i\in\mathcal{B}_{0}$ consist of the critical load factors $\lambda_{i}$ and the associated buckling modes $\boldsymbol{\varphi}_{i}$, normalized such that $\boldsymbol{\varphi}^{T}_{j}\mathbf{K}_{\sigma}\boldsymbol{\varphi}_{i} = -\delta_{ji}$.

In this work we are particularly concerned with the fundamental buckling load factor $\lambda_{1}$, associated with the critical load $\mathbf{f}_{\rm cr} = \lambda_{1}\mathbf{f}_{0}$. Therefore, from a computational point of view it is more convenient to refer to the eigenvalue equation
\[
  \left[ \mathbf{K}_{\sigma}\left( \mathbf{u}_{0} \right) - \mu\mathbf{K} \right] \boldsymbol{\varphi} = \mathbf{0} \: ,
  \qquad \boldsymbol{\varphi} \neq \mathbf{0}
\]
which is equivalent to \eqref{eq:GEPforMu} given the relationship $\lambda = -1/\mu$, such that $\lambda_{1}$ is associated with the minimum algebraic value of $\mu$, say $\mu_{1}$.

For density--based topology optimization \citep{book:bendsoe-sigmund_2004} performed on a regular grid, we consider the design variables $\mathbf{x} = \left\{ x_{e} \right\}^{m}_{e=1}$, belonging to the set
\begin{equation}
 \label{eq:FeasibleSet}
 \mathcal{F} := \left\{ \mathbf{x} \in \left[ 0, 1 \right]^{m} \, , \quad \frac{1}{m f }\sum^{m}_{e=1} x_{e} - 1 \leq 0 \right\}
\end{equation}
where $ f $ is the allowed volume fraction, and the design--dependent matrices $\mathbf{K}\left( \mathbf{x} \right)$ and $\mathbf{K}_{\sigma}\left( \mathbf{x}, \mathbf{u}_{0} \right)$ are assembled from the element ones.

The latter are parametrized by interpolation functions as $\mathbf{k}_{e} = h_{1}\left( x_{e} \right)\mathbf{k}_{0}$ and $\mathbf{g}_{e}\left( \mathbf{x}, \mathbf{u} \right) = h_{2}\left( x_{e} \right)\mathbf{g}_{0}\left( \mathbf{u}_{0e} \right)$, where $\mathbf{u}_{0e}$ denotes the restriction of the global displacement vector to the element level. In the following we consider
\begin{equation}
 \label{eq:InterpolationFunctions}
  \begin{aligned}
   h_{1}\left( x_{e} \right) & = E_{0} + x^{p}_{e}\left( E_{1} - E_{0} \right) \\
   h_{2}\left( x_{e} \right) & = x^{p}_{e}E_{1}
  \end{aligned}
\end{equation}
where $E_{1}$ is the Young modulus of the solid, $E_{0}$ that of the void and $p$ is the penalization factor.

\autoref{eq:InterpolationFunctions} has proven to be an effective choice against artificial buckling modes \citep{gao-ma_15a,thomsen-etal_18a}, at least for the material contrast of interest ($E_{1}/E_{0} = 10^{6}$).

The problem of minimizing the linear compliance with an imposed lower bound $\overline{P_{c}}$ on the fundamental buckling load factor reads
\begin{equation}
 \label{eq:OP-MinComplianceWithBucklingConstr}
  \begin{cases}
              & \min\limits_{\mathbf{x}\in\mathcal{F}} J := \mathbf{u}^{T}_{0}\mathbf{f} = \mathbf{u}^{T}_{0}\mathbf{K}\mathbf{u}_{0}\\
   {\rm s.t.} & \min\limits_{i\in\mathcal{B}}\lambda_{i} \geq \overline{P_{c}}
  \end{cases}
\end{equation}
where $\mathcal{B} \subset \mathcal{B}_{0}$ is the subset of eigenvalues considered, which should be large enough to consider all the relevant buckling modes and to produce a smooth behavior of the optimizer \citep{bruyneel-etal_08a,dunning-etal_16a}. The side constraint in \eqref{eq:OP-MinComplianceWithBucklingConstr} can be replaced by the set
\begin{equation}
 \label{eq:EigConstraintMU}
 \overline{P_{c}}\alpha^{(1-i)}\mu_{i} + 1
 \geq 0 \ ,
 \qquad i\in\mathcal{B}
\end{equation}
where the number $\alpha \geq 1$ serves to introduce small gaps between the eigenvalues, preventing them from completely coalescing \citep{book:bendsoe-sigmund_2004}.

The solution of \eqref{eq:OP-MinComplianceWithBucklingConstr} with an iterative, two--level, gradient based approach \citep{book:haftka2012}, requires the solution of \eqref{eq:LinearizedEquilibriumEquation} and \eqref{eq:GEPforMu} and the computation of sensitivities used for updating the design variables. For a simple eigenvalue $\lambda_{i}$, the sensitivity with respect to the design variable $x_{e}$ reads \citep{rodrigues-etal_95a}
\begin{equation}
 \label{eq:SensitivityLambda}
  \frac{\partial \lambda_{i}}{\partial x_{e}} = \boldsymbol{\varphi}^{T}_{i} \left( \frac{\partial \mathbf{K}}{\partial x_{e}} + \lambda_{i} \frac{\partial \mathbf{K}_{\sigma}}{\partial x_{e}} \right)
  \boldsymbol{\varphi}_{i} - \lambda_{i}\mathbf{v}^{T}\frac{\partial \mathbf{K}}{\partial x_{e}}\mathbf{u}_{0}
\end{equation}
where $\mathbf{v}$ is obtained by solving the adjoint system
\begin{equation}
 \label{eq:AdjointSystem}
 \mathbf{K}\mathbf{v} = \boldsymbol{\varphi}^{T}_{i}\left[ \nabla_{\mathbf{u}}\mathbf{K}_{\sigma} \right]\boldsymbol{\varphi}_{i}
\end{equation}

The first term in \eqref{eq:SensitivityLambda} is formally the same as the sensitivity for a dynamic eigenvalue, and in the following we will refer to it as the ``frequency--like'' term. The second term is the adjoint term, accounting for the dependence of the stress stiffness matrix on the stress level in the prebuckling solution, and the variation of this as the design is changed \citep{rodrigues-etal_95a}.

It is sometimes claimed that the latter is of minor importance and can be neglected \citep{mateus-etal_97a,munk-etal_16a}. However, ignoring this term results in the use of inconsistent gradients and may potentially lead to incorrect designs, as will be discussed later.

\autoref{eq:SensitivityLambda} does not give an univocal value for a repeated eigenvalue, as more than one eigenvector becomes associated with it. In this case the definition of the sensitivity becomes more involved, calling for the use of subgradients, and we refer to \citet{rodrigues-etal_95a} for details. In this work we do not consider such more difficult treatment, either because setting $\alpha > 1$ in \eqref{eq:EigConstraintMU} will ensure a certain separation of the eigenvalues or, as an alternative, we will make use of the aggregation functions discussed below.

\begin{figure*}[t]
  \begin{minipage}[b]{0.6\linewidth}
    \centering
    \includegraphics[scale = 0.5, keepaspectratio]
         {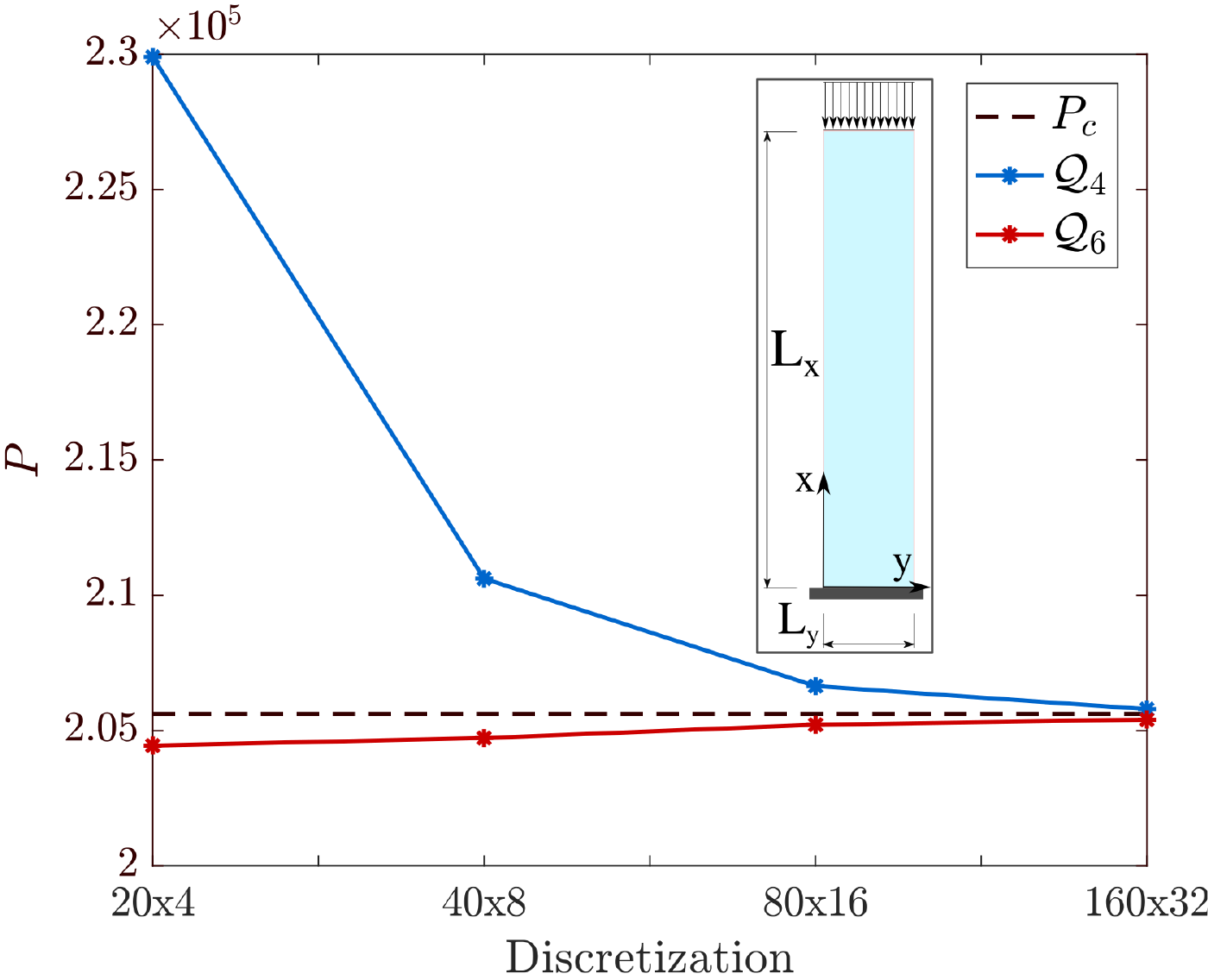}
  \end{minipage}%
  \begin{minipage}[b]{0.4\linewidth}
    \centering
    \begin{tabular}[b]{ccc}
     \hline\noalign{\smallskip}
     $m_{x} \times m_{y}$ & \multicolumn{2}{c}{$e_{r}$} \\
     \noalign{\smallskip}\hline
      & $\mathcal{Q}_{4}$ & $\mathcal{Q}_{6}$ / PS \\
     \noalign{\smallskip}\hline
     $ 10 \times  2$ &            $0.494$ & $-6.9 \cdot 10^{-3}$ \\
     $ 10 \times  4$ &            $0.493$ & $-6.9 \cdot 10^{-3}$ \\
     $ 20 \times  2$ &            $0.119$ & $-6.6 \cdot 10^{-3}$ \\
     $ 20 \times  4$ &            $0.118$ & $-5.7 \cdot 10^{-3}$ \\
     $ 40 \times  8$ &            $0.024$ & $-5.1 \cdot 10^{-3}$ \\
     $ 80 \times 16$ & $5.0\cdot 10^{-3}$ & $-2.6 \cdot 10^{-3}$ \\
     $160 \times 32$ & $8.6\cdot 10^{-4}$ & $-2.1 \cdot 10^{-3}$ \\
     \noalign{\smallskip}\hline
    \end{tabular}
\end{minipage}
\caption{Accuracy of $\lambda_{1}$ approximations obtained by using $\mathcal{Q}_{4}$ or $\mathcal{Q}_{6}$ discretizations for the compressed column example. Results obtained with the $\mathcal{Q}_{6}$ element are identical (up to numerical precision) to those obtained with the Pian--Sumihara (PS) element.}
\label{fig:AccuracyCompressedColumn}
\end{figure*}

For completeness, we also report the sensitivity for the eigenvalues $\mu_{i}$
\begin{equation}
 \label{eq:SensitivityMu}
  \frac{\partial \mu_{i}}{\partial x_{e}} = \boldsymbol{\varphi}^{T}_{i} \left( \frac{\partial \mathbf{K}_{\sigma}}{\partial x_{e}} - \mu_{i} \frac{\partial \mathbf{K}}{\partial x_{e}} \right)
  \boldsymbol{\varphi}_{i} - \mathbf{v}^{T}\frac{\partial \mathbf{K}}{\partial x_{e}}\mathbf{u}_{0}
\end{equation}
which is the one actually adopted when referring to \eqref{eq:EigConstraintMU}. The adjoint vector $\mathbf{v}$ is again obtained by \eqref{eq:AdjointSystem}.

\subsection{The Finite Element Approximation}
 \label{sSec:ShortCommentFEA}

Low order elements, often preferred for density--based TO, may result in inaccurate representation of stresses, which are of primary importance in buckling analysis. In this regard, conforming elements are known to suffer from shear locking and therefore to poorly represent a bending dominated behavior. Although this is usually neglected in compliance design, since the optimized topology mainly consists of tension and compression members, it might play a role for buckling.

A number of accurate low order elements have been developed in the past, mainly relying on the enrichment of the compatible displacement field with Incompatible Modes \citep{turner-etal_56a,wilson-etal_73a}, or derived from mixed formulations, as Hybrid Stress models \citep{pian_64a,pian-sumihara_84a} or Enhanced Assumed Strains models \citep{simo-rifai_90a}. The main advantage of all these models is to substantially increase the accuracy of the stress approximation on coarse meshes without introducing new nodes or explicit DOFs.

The classic hybrid stress element of \citet{pian-sumihara_84a} has been applied to buckling TO by \citet{gao-ma_15a,gao-etal_17a}, however without discussing its advantages over the conforming $\mathcal{Q}_{4}$ element. Here we shortly analyze the behavior of a popular incompatible element, the so--called $\mathcal{Q}_{6}$ element, proposed by \cite{wilson-etal_73a}. We remark that the conclusions we draw also apply to the Pian element, since the two have proven to be equivalent by \citet{froier-etal_74a}.

\begin{figure*}[t]
 \centering
  \subfloat[]{
   \includegraphics[scale = 0.225, keepaspectratio]
   {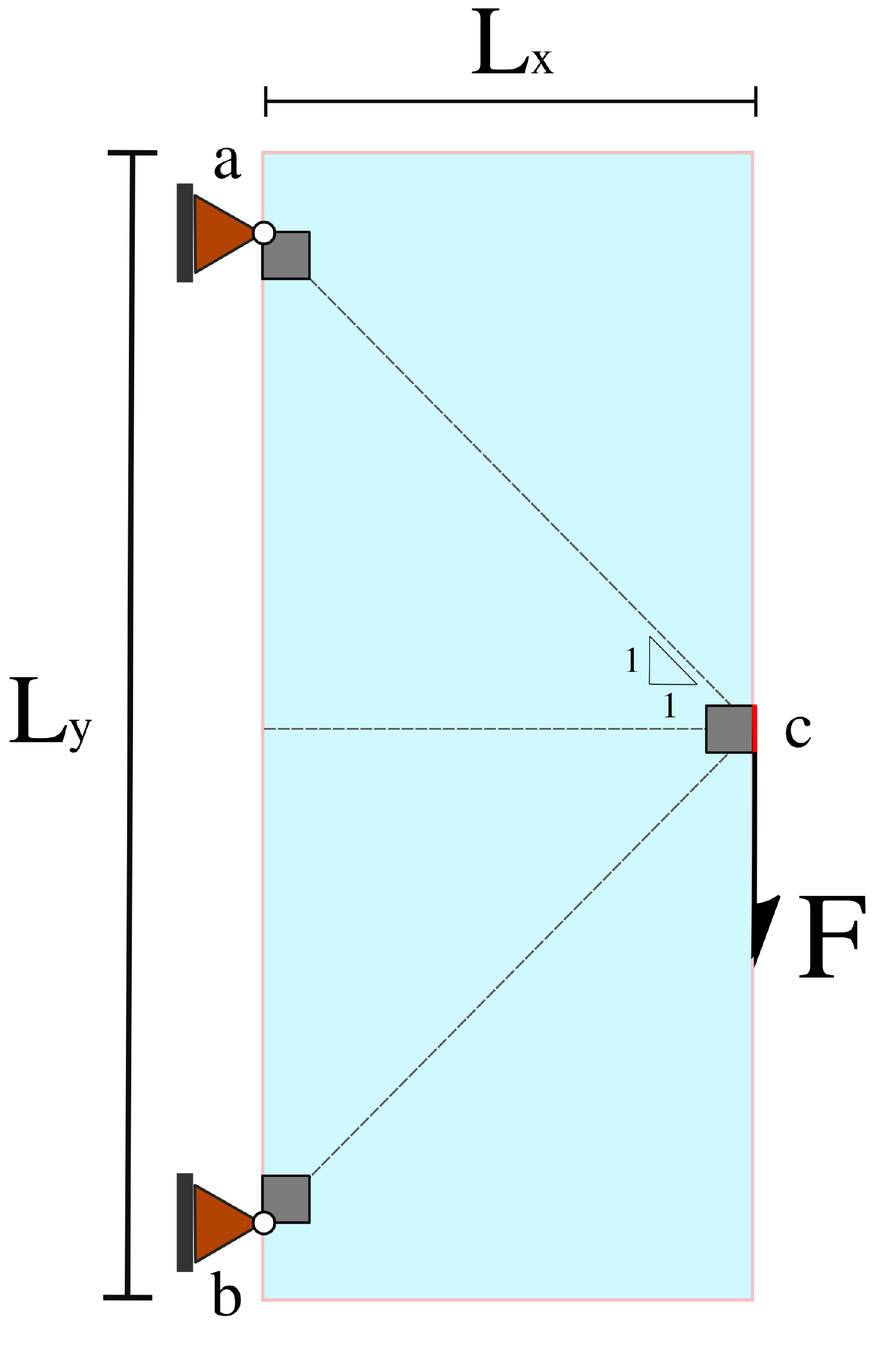}} \quad
  \subfloat[$\overline{P_{c}} = 0$]{
   \includegraphics[scale = 0.225, keepaspectratio]
   {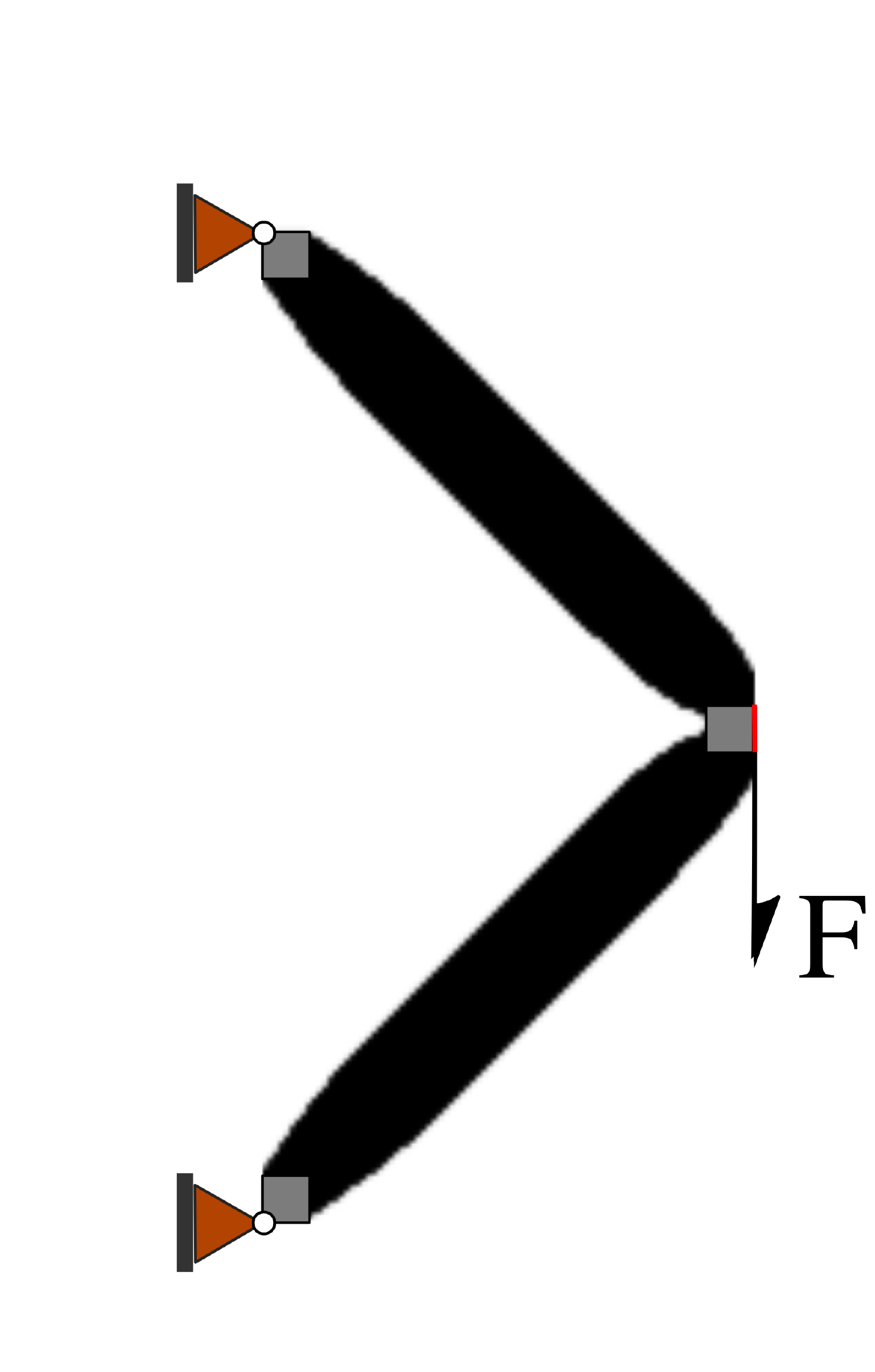}} \quad
  \subfloat[$\log\left( \phi_{e}/\phi_{\rm max} \right)$]{
   \includegraphics[scale = 0.225, keepaspectratio]
   {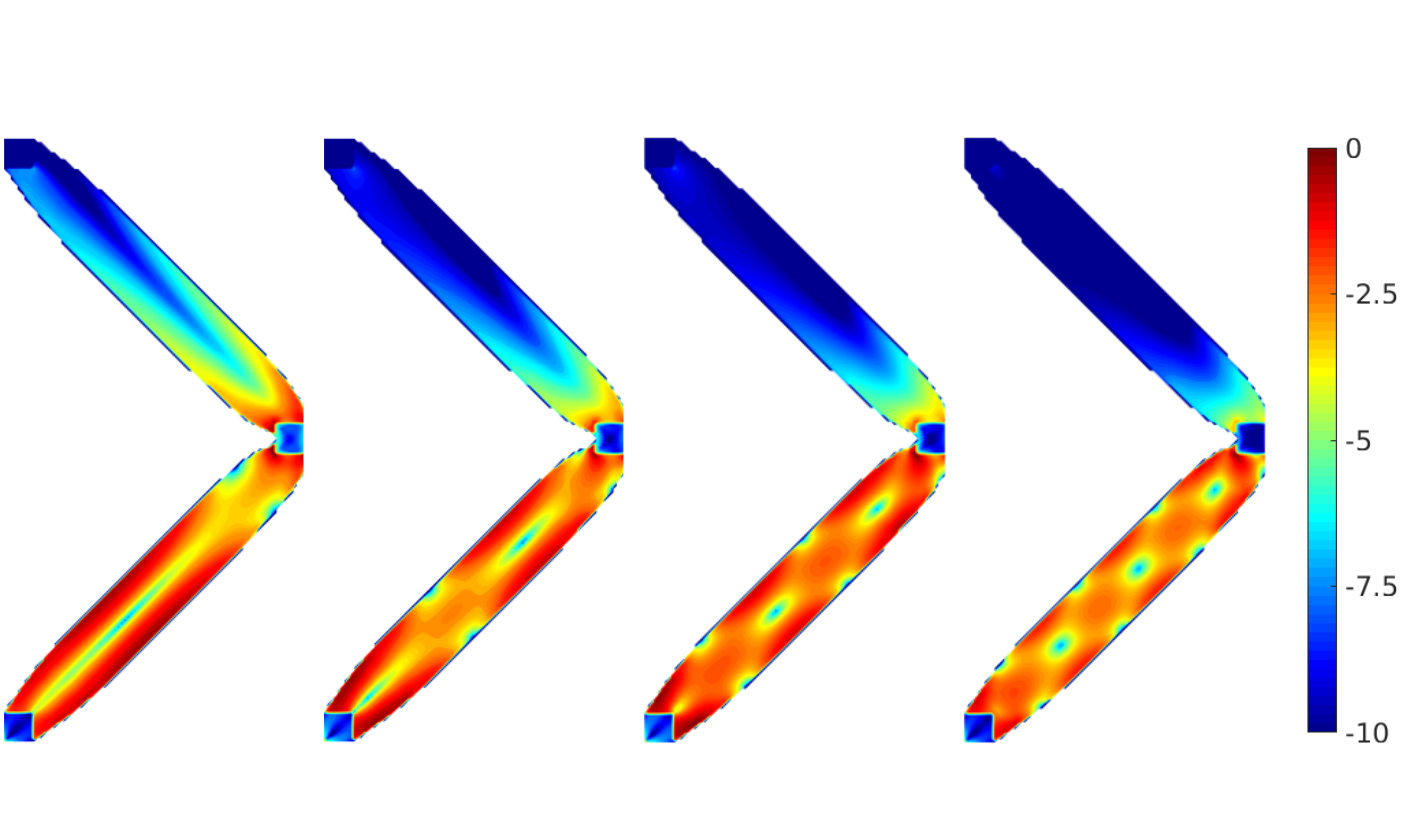}}
 \caption{Geometrical setting for the considered example (a) compliance design (b) and distribution of the strain energy density for the lowest four buckling modes (c)}
 \label{fig:GeometricalSettingAndComplianceDesing}
\end{figure*}

Let us consider the following example, recently considered by \citet{pedersen-pedersen_18b}: a cantilever beam with dimensions $L_{x} = 10$, $L_{y} = 1$ and thickness $t = 0.005$ (all measures non-dimensioanalized), subjected to uniform compression at the tip. The elastic modulus is $E = 2\cdot 10^{11}$ and we set $\nu = 0$. Due to this parameter choice, the 2D plane stress model is equivalent to the 1D beam model, and the fundamental buckling load can be computed by the closed form expression

\begin{equation}
  P_{c} = \frac{E I_{y}}{4 L^{2}_{x}} \approx 2.05617 \cdot 10^{6}
\end{equation}
where $I_{y} = t L^{3}_{y}/12$ is the modulus of inertia for in--plane bending (see the sketch in \autoref{fig:AccuracyCompressedColumn}). We choose a reference load vector with magnitude $\left\| \mathbf{f}_{0} \right\| = 2.5 \cdot 10^{5}$. The convergence behavior and relative errors associated with the $\mathcal{Q}_{4}$ and the $\mathcal{Q}_{6}$ element discretizations are reported in  \autoref{fig:AccuracyCompressedColumn}.

The use of $\mathcal{Q}_{4}$ conforming elements results in considerable errors on coarse discretizations, with an overestimation of $P_{c}$ by almost $50\%$ for a column described by $10 \times 4$ elements. The error is reduced to the order of $10^{-2}$ and $10^{-3}$ when increasing the number of elements by 8 and 32 times, respectively. 

On the other hand, using the $\mathcal{Q}_{6}$ element sensibly improves the accuracy on coarse discretizations. Moreover, it is important to observe that now the convergence is attained from below, since with the incompatible finite element model, the discretized structure is softer than the real one. For the $\mathcal{Q}_{6}$ element we acknowledge a very mild convergence as the mesh is refined and therefore the performance of the two elements becomes very similar on fine meshes.

From the above discussion it should be clear that buckling modes might be inaccurately predicted for optimized designs consisting of thin slender bars, described by few $\mathcal{Q}_{4}$ elements. The performance of the two elements for the analysis of optimized designs and their influence on the optimization itself will be further discussed in \autoref{sSec:ComparisonResultsQ6}.

\subsection{Formulation by aggregation functions}
 \label{sSec:FormulationWithAggregationFunctions}

We can replace \eqref{eq:EigConstraintMU} with a single constraint
\begin{equation}
 \label{eq:SingleConstraint}
  \overline{P_{c}}M\left[ \mu_{i} \right] + 1 \geq 0
\end{equation}
where $M\left[ \cdot \right]$ is a suitable aggregation function \citep{chen-etal_04a}. Various functions can be adopted in this regard, and the most popular ones in the context of TO are the $\rho$--norm
\begin{equation}
 \label{eq:pNormAggregation}
  M_{\rho}\left[ \mu_{i} \right] =
  \left( \sum_{i\in\mathcal{B}}\left( - \mu_{i} \right)^{\rho} \right)^{1/\rho}
\end{equation}
where we used that the $\mu_{i}$ of interest are strictly negative, and the KS function \citep{kreisselmeier-steinhauser_79a}
\begin{equation}
 \label{eq:ksAggregation}
   M_{KS}\left[ \mu_{i} \right] = \mu_{1} + \frac{1}{\rho} \ln\left( \sum_{i\in\mathcal{B}}
   e^{\rho\left( \mu_{i} - \mu_{1} \right)} \right)
\end{equation}

Both \eqref{eq:pNormAggregation} and \eqref{eq:ksAggregation} produce an upper bound, smooth approximation to $\mu_{1} = \max\limits_{i\in\mathcal{B}}|\mu_{i}|$, where the degree of smoothness is governed by the aggregation parameter $\rho$. Therefore, we indirectly obtain a smooth lower bound to $\lambda_{1}$. Ideally, $\mu_{1}$ is recovered as $\rho\rightarrow\infty$.

The KS function is usually preferred if a high $\rho$ value has to be used, due to its superior numerical stability. For a more complete discussion about the properties of these function and for an overview on some more general aggregation strategies we refer, for example, to \citet{raspanti-etal_00a,kennedy-hicken_15a}.

The sensitivity of the aggregated constraint \eqref{eq:SingleConstraint} can be easily obtained as
\begin{equation}
 \label{eq:pNormAggregationSens}
  \frac{\partial M_{\rho}}{\partial x_{e}} =
  - M^{(1 - \rho)}_{\rho}
  \sum_{i\in\mathcal{B}} \mu_{i}^{(\rho-1)}
  \frac{\partial \mu_{i}}{\partial x_{e}}
\end{equation} %
when the $\rho$--norm is used and
\begin{equation}
 \label{eq:ksAggregationSens}
  \frac{\partial M_{{\rm KS}}}{\partial x_{e}} =
  \frac{1}
  {\sum_{i\in\mathcal{B}} e^{\rho \mu_{i}}}
  \sum_{i\in\mathcal{B}} e^{\rho \mu_{i}}
  \frac{\partial \mu_{i}}{\partial x_{e}}
\end{equation}
when using the KS function. In the above the sensitivities $\partial\mu_{i}/\partial x_{e}$ are given by \eqref{eq:SensitivityMu}.

One of the advantages of using the aggregation functions is the uniqueness of the gradient defined through \eqref{eq:pNormAggregationSens} (viz. \eqref{eq:ksAggregationSens}), even when eigenvalues are repeated. A proof of this, based on the argument of \citet{gravesen-etal_11a}, can be found in \citet{torii-faria_17a}, for the $\rho$--norm function. The same argument trivially extends also to the KS function, due to the continuity of the exponential and logarithm functions involved.

\begin{figure*}[tb]
 \centering
  \subfloat[$\overline{P_{c}} = 0.5$]{
   \includegraphics[scale = 0.2, keepaspectratio]
   {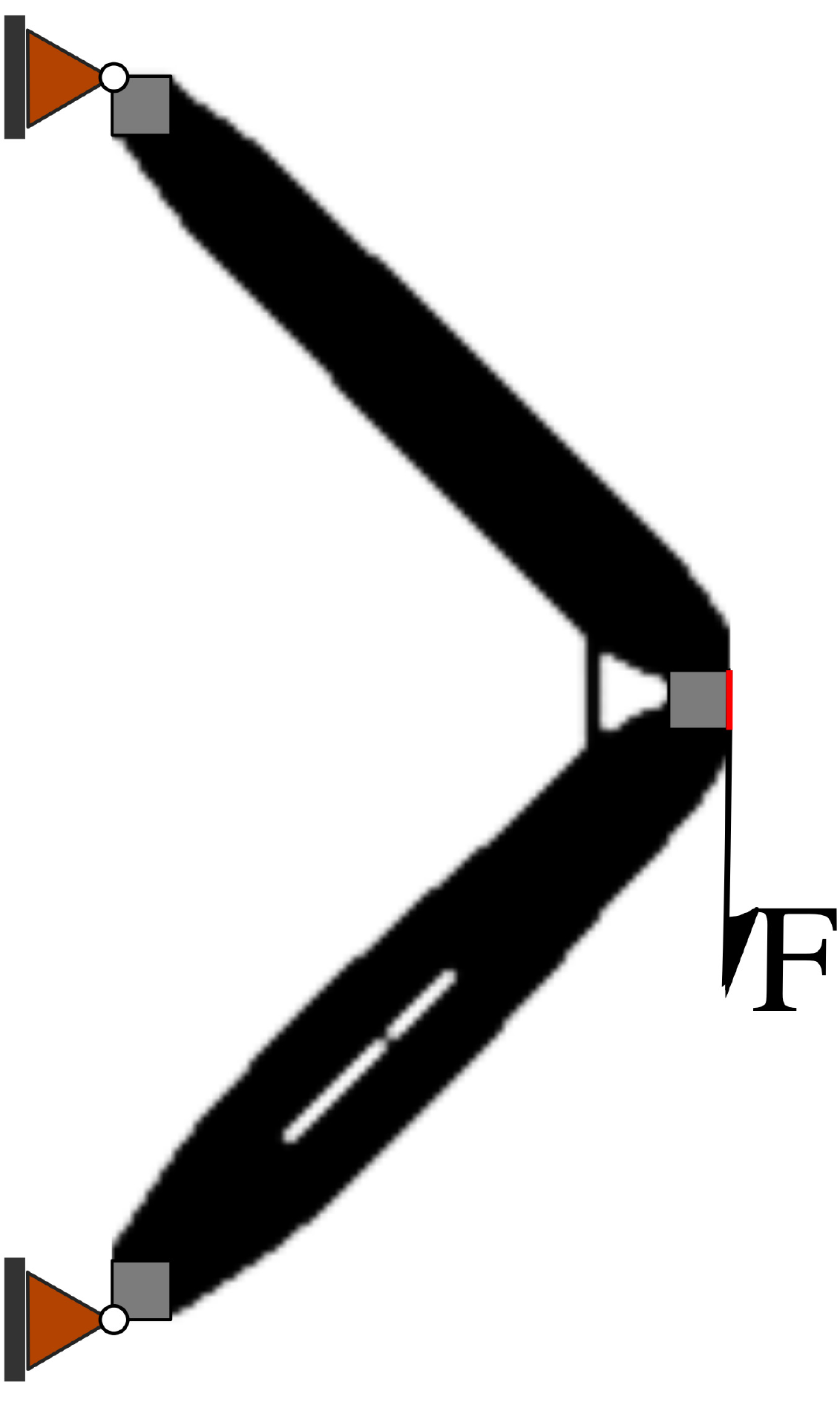}} \quad
  \subfloat[$\overline{P_{c}} = 0.75$]{
   \includegraphics[scale = 0.2, keepaspectratio]
   {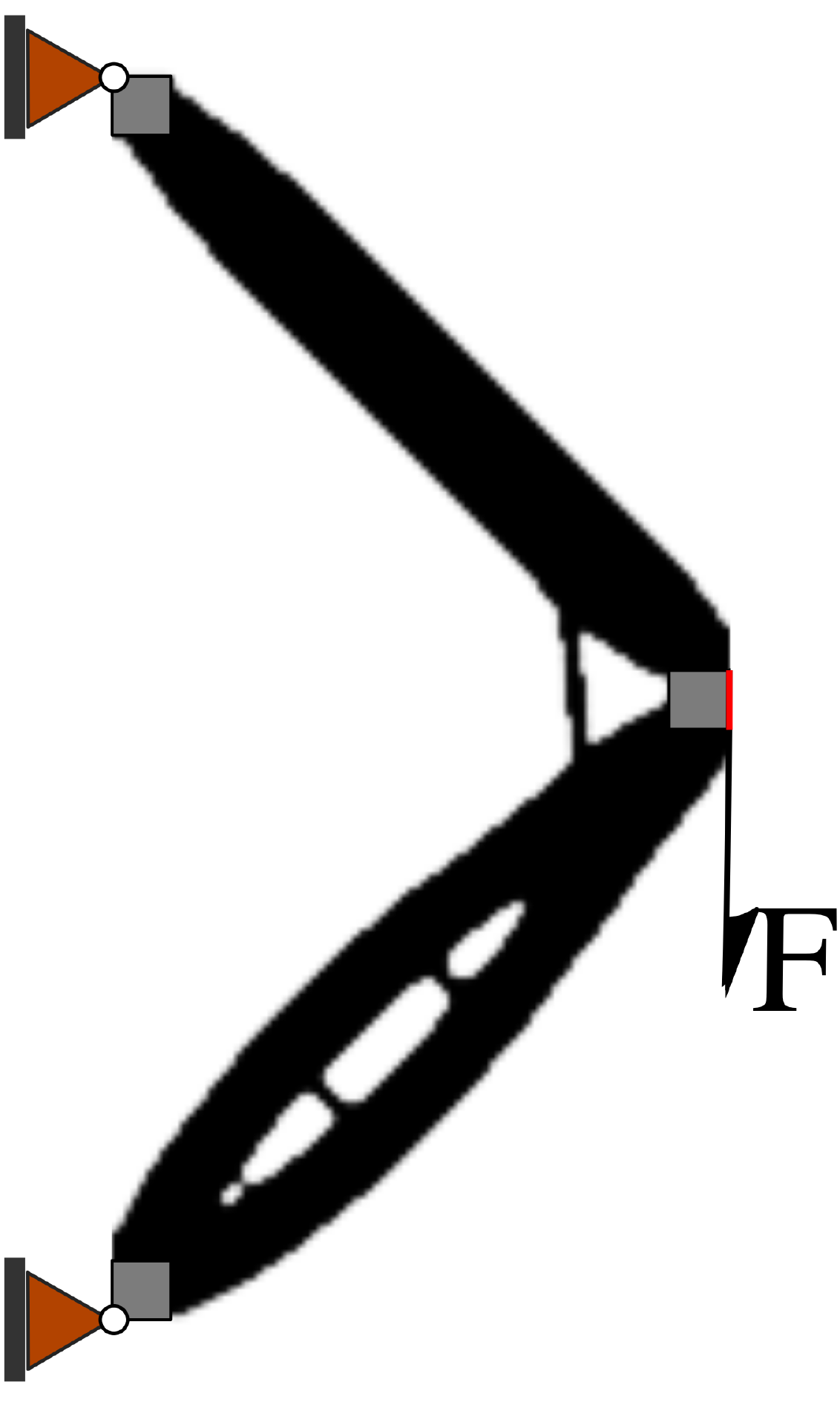}} \quad
  \subfloat[$\overline{P_{c}} = 1.0$]{
   \includegraphics[scale = 0.2, keepaspectratio]
   {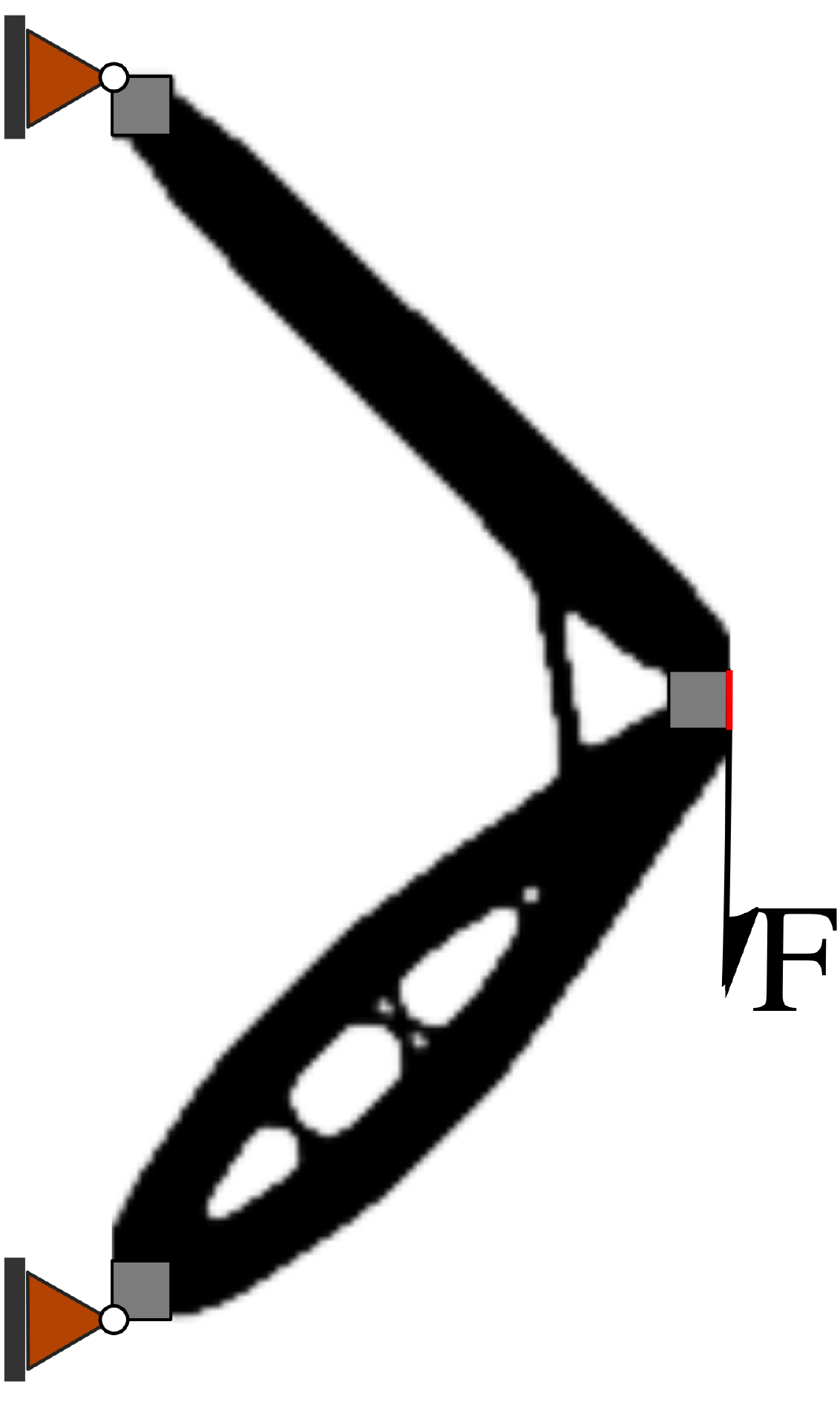}} \quad
  \subfloat[$\overline{P_{c}} = 1.25$]{
   \includegraphics[scale = 0.2, keepaspectratio]
   {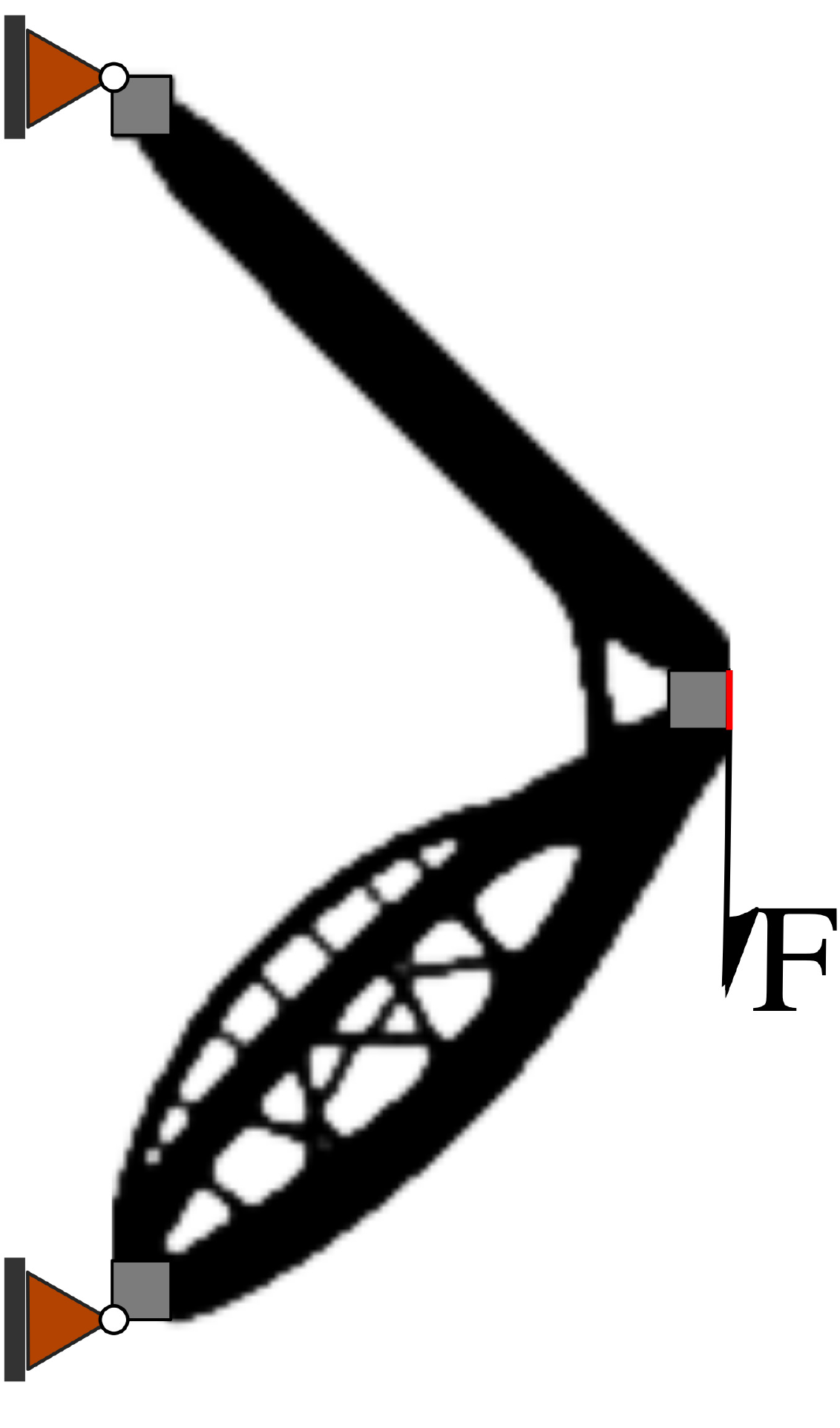}} \quad
  \subfloat[$\overline{P_{c}} = 1.5$]{
   \includegraphics[scale = 0.2, keepaspectratio]
   {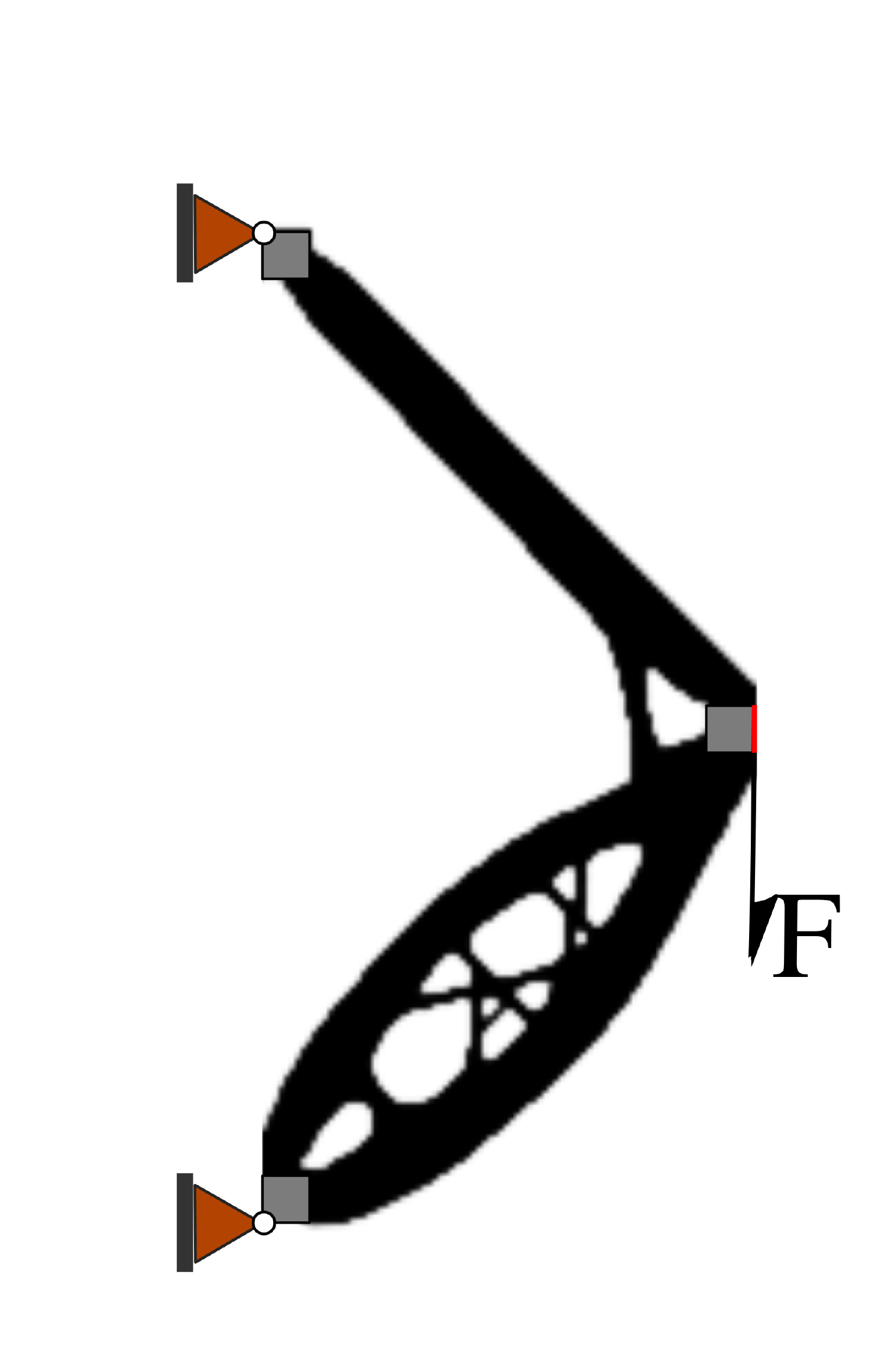}} \\
  \subfloat[$\overline{P_{c}} = 1.75$]{
   \includegraphics[scale = 0.2, keepaspectratio]
   {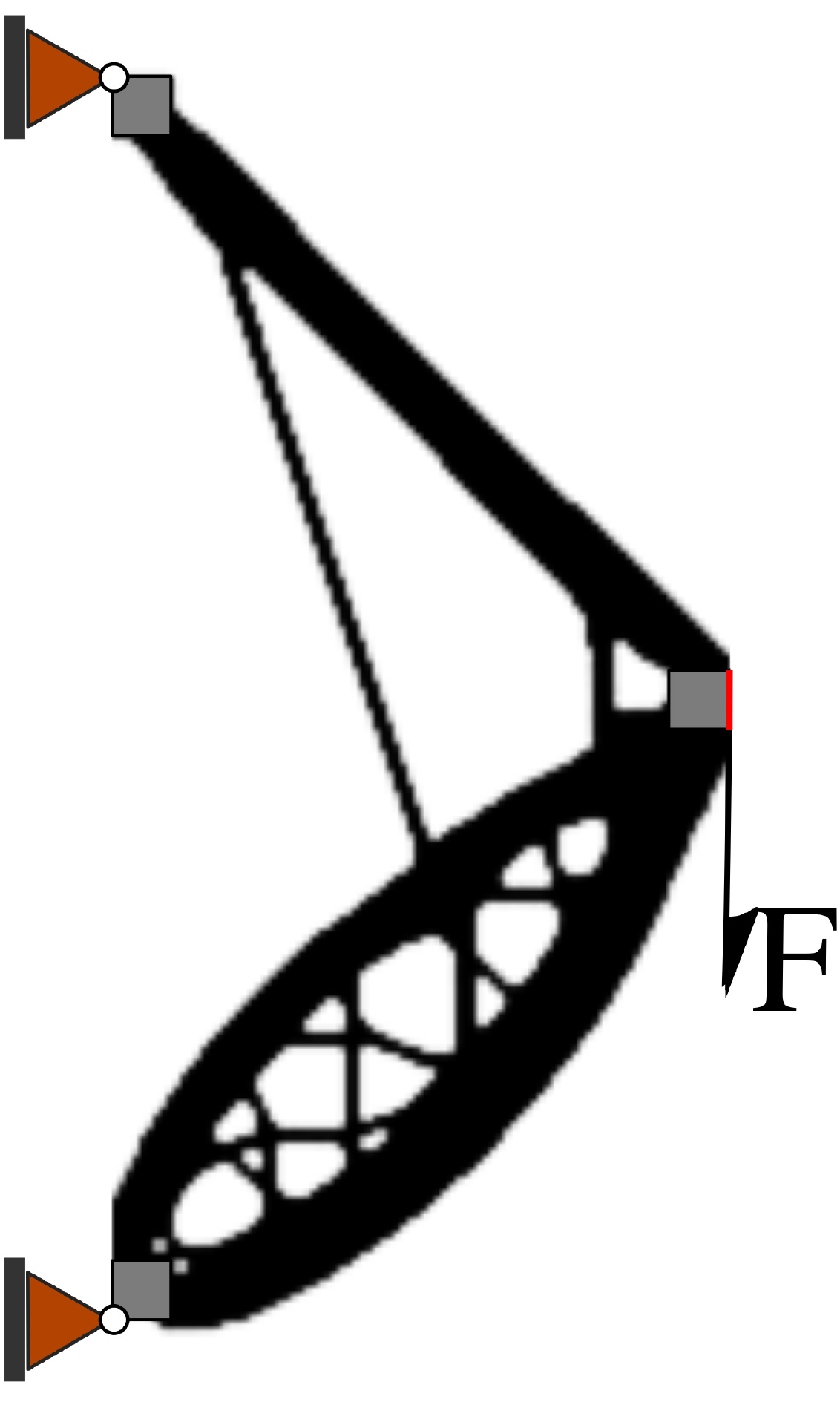}} \quad
  \subfloat[$\overline{P_{c}} = 1.8$]{
   \includegraphics[scale = 0.2, keepaspectratio]
   {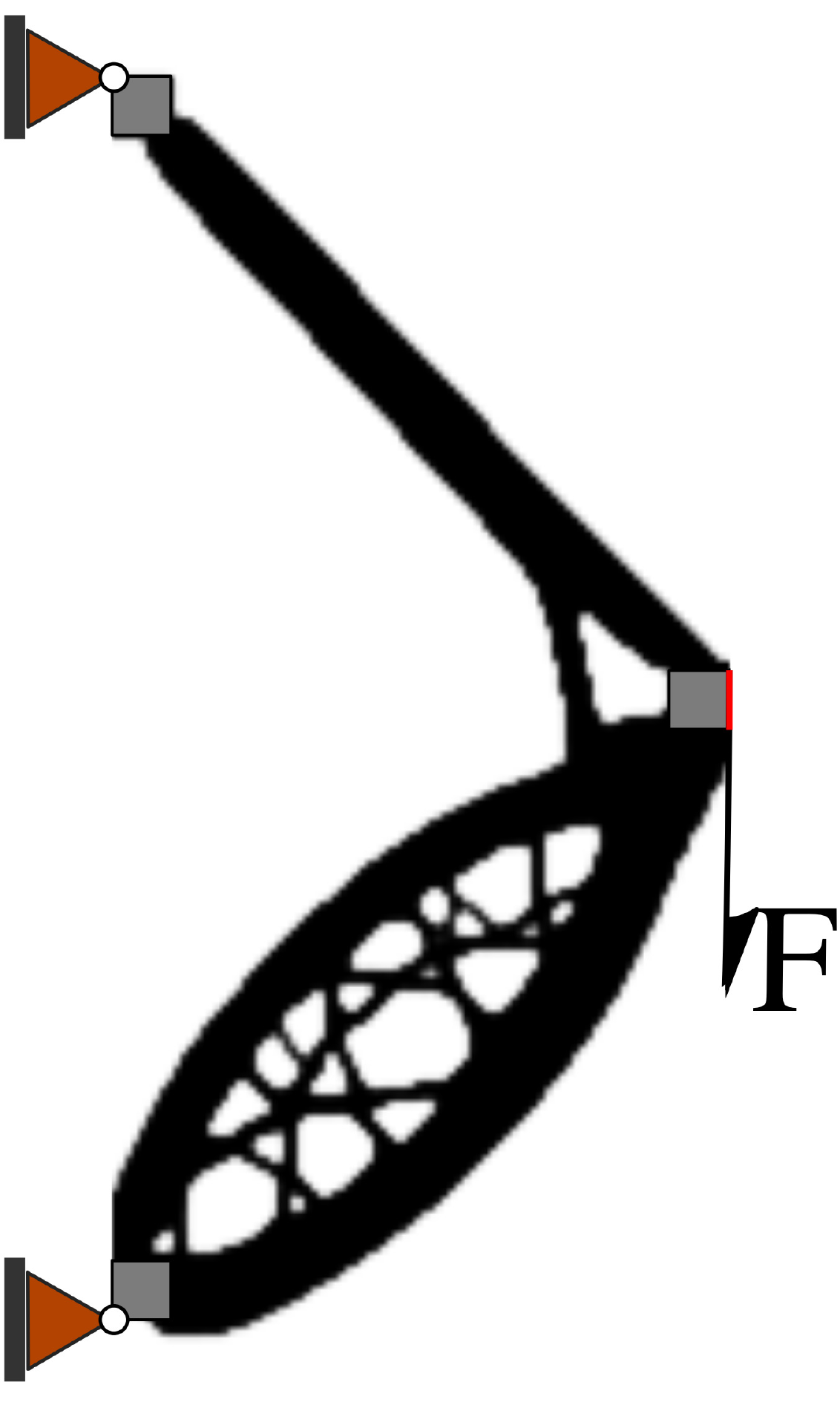}} \quad
  \subfloat[$\overline{P_{c}} = 2.0$]{
   \includegraphics[scale = 0.2, keepaspectratio]
   {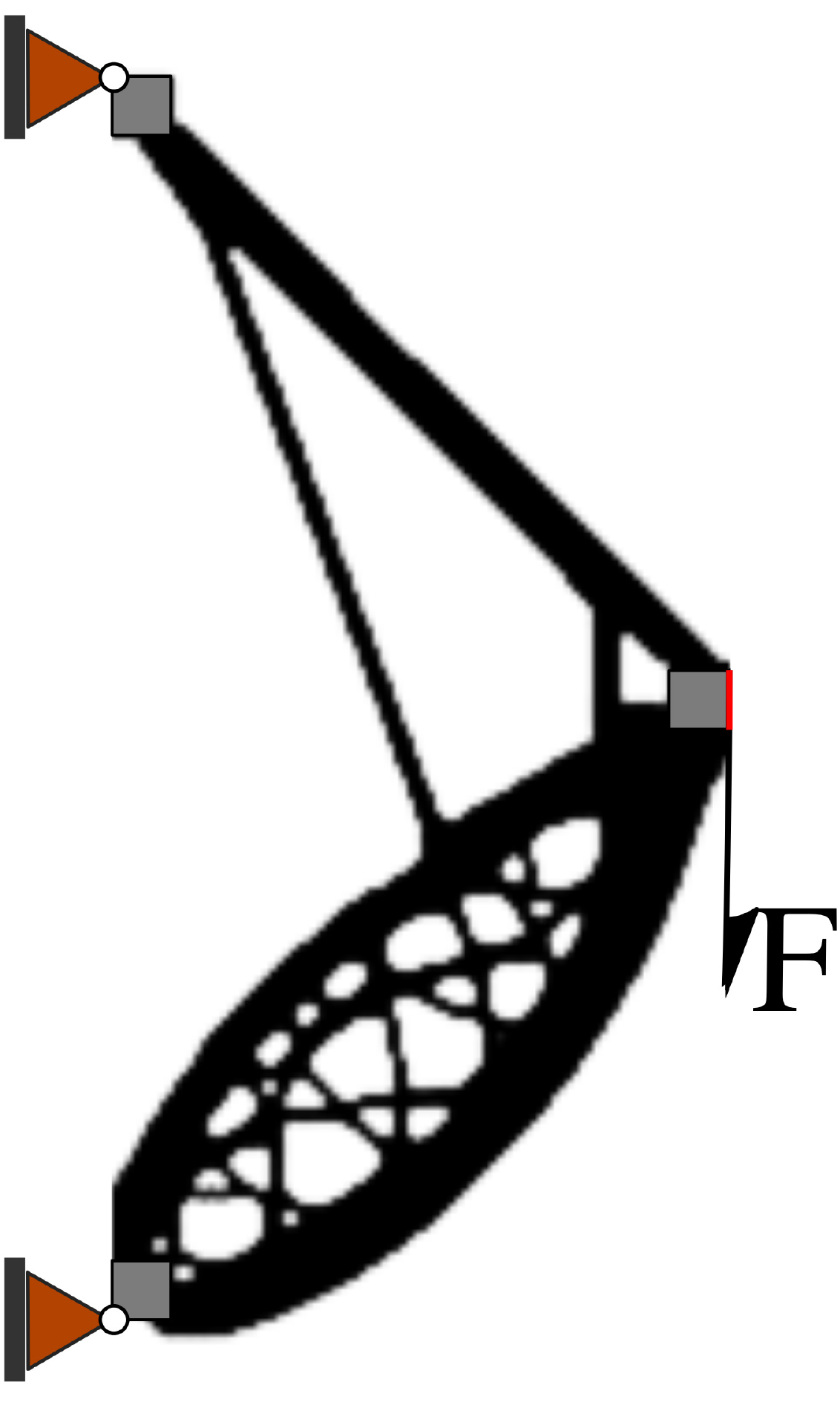}} \quad
  \subfloat[$\overline{P_{c}} = 2.0^{\ast}$]{
   \includegraphics[scale = 0.2, keepaspectratio]
   {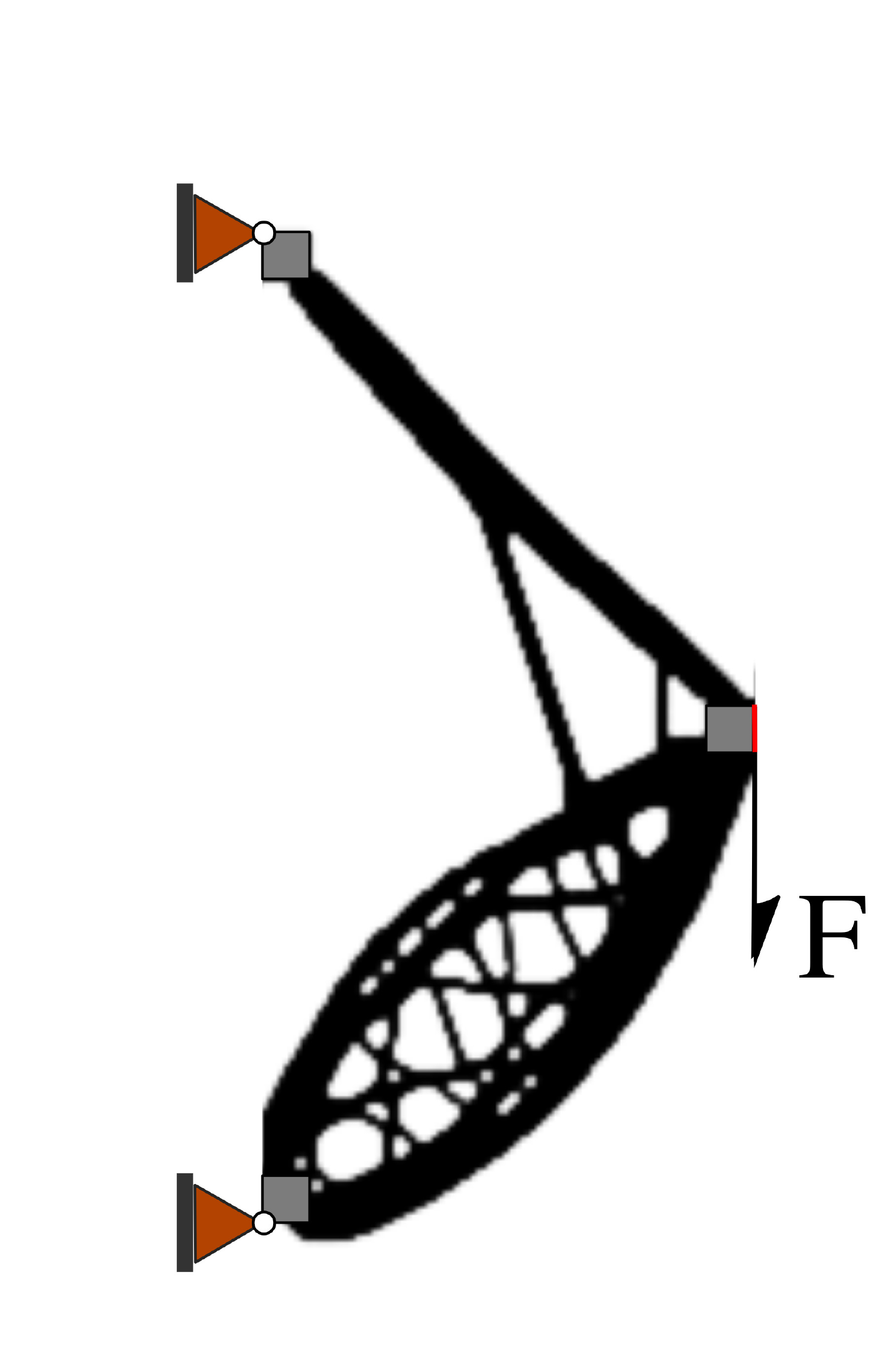}} \quad
  \subfloat[$\overline{P_{c}} = 2.0^{\ast\ast}$]{
   \includegraphics[scale = 0.2, keepaspectratio]
   {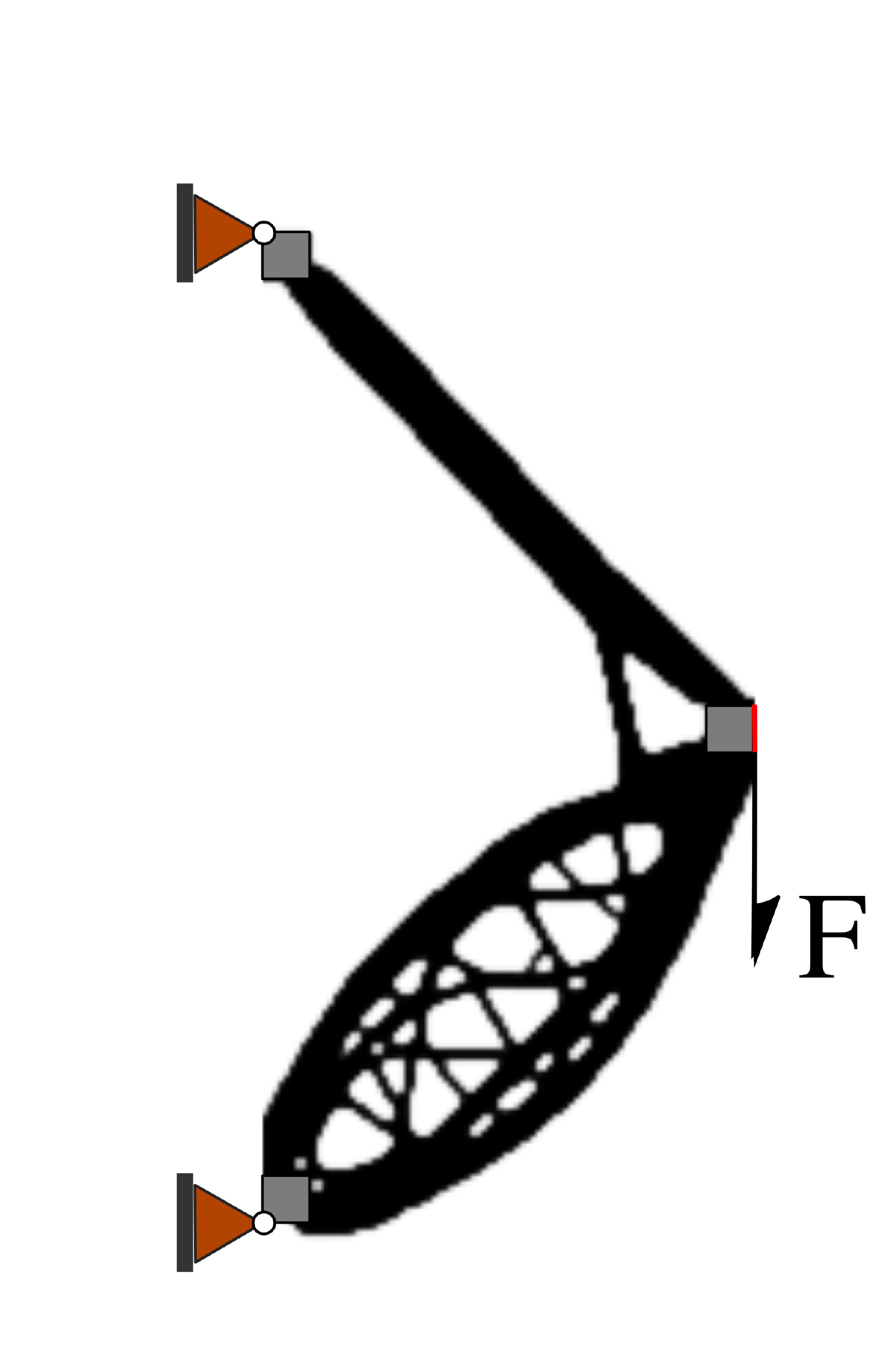}}
 \caption{Optimized designs obtained for the example described in \autoref{fig:GeometricalSettingAndComplianceDesing} (a) and for increasing values of the lower bound on the fundamental buckling load factor. Designs marked with $(^{\ast})$ and $(^{\ast\ast})$ have been obtained starting from non--uniform initial guesses: (i) is obtained from the compliance design of \autoref{fig:OptimizedDesignsIncreasingPc} (b) and (j) from the design corresponding to $\overline{P_{c}} = 1$ (see \autoref{fig:OptimizedDesignsIncreasingPc} (c))}
 \label{fig:OptimizedDesignsIncreasingPc}
\end{figure*}

\section{Numerical example}
 \label{Sec:NumExp}

We refer to the configuration of \autoref{fig:GeometricalSettingAndComplianceDesing} (a), where the gray regions are fixed to be solid in the optimization. The load, having total magnitude $F = 0.02$, is distributed over the edge of the gray region on the right side, centered in the mid point $c$, and the pinned points $a$ and $b$ on the left side are on the $45^{\circ}$ line from $c$. We here assume the load to be unidirectional and pointing downward, and therefore a single load case needs to be considered and only positive eigenvalues of \eqref{eq:GEPforMu} are accounted for.

The discretization is set up with ${\rm \Omega}_{1} = 90 \times 210$ elements and the solid regions consist of $9\times 9$ elements. We consider $E_{1} = 1$ and $E_{0} = 10^{-6}$ in \eqref{eq:InterpolationFunctions}, and the Poisson ratio is $\nu = 0.3$. The penalization parameter is initialized with $p = 1$ and then raised in steps of $\Delta p = 0.25$ each 25 optimization steps, up to the value $p = 6$. With this continuation approach we do not start with a too weak design, and we increase the chances of ending in a good local minimum. More advanced, adaptive continuation strategies can be used and for this we refer, for example, to \citet{labanda-stolpe_15a,gao-etal_17a}. A projection with threshold $\eta = 0.5$ and sharpness parameter $\beta = 6$, based on a filter radius of $r_{\rm min} = 2$ elements \citep{sigmund_07a,wang-etal_11a}, is applied to control the structural features size.

\begin{figure}[t]
 \centering
  \includegraphics[scale = 0.525, keepaspectratio]
   {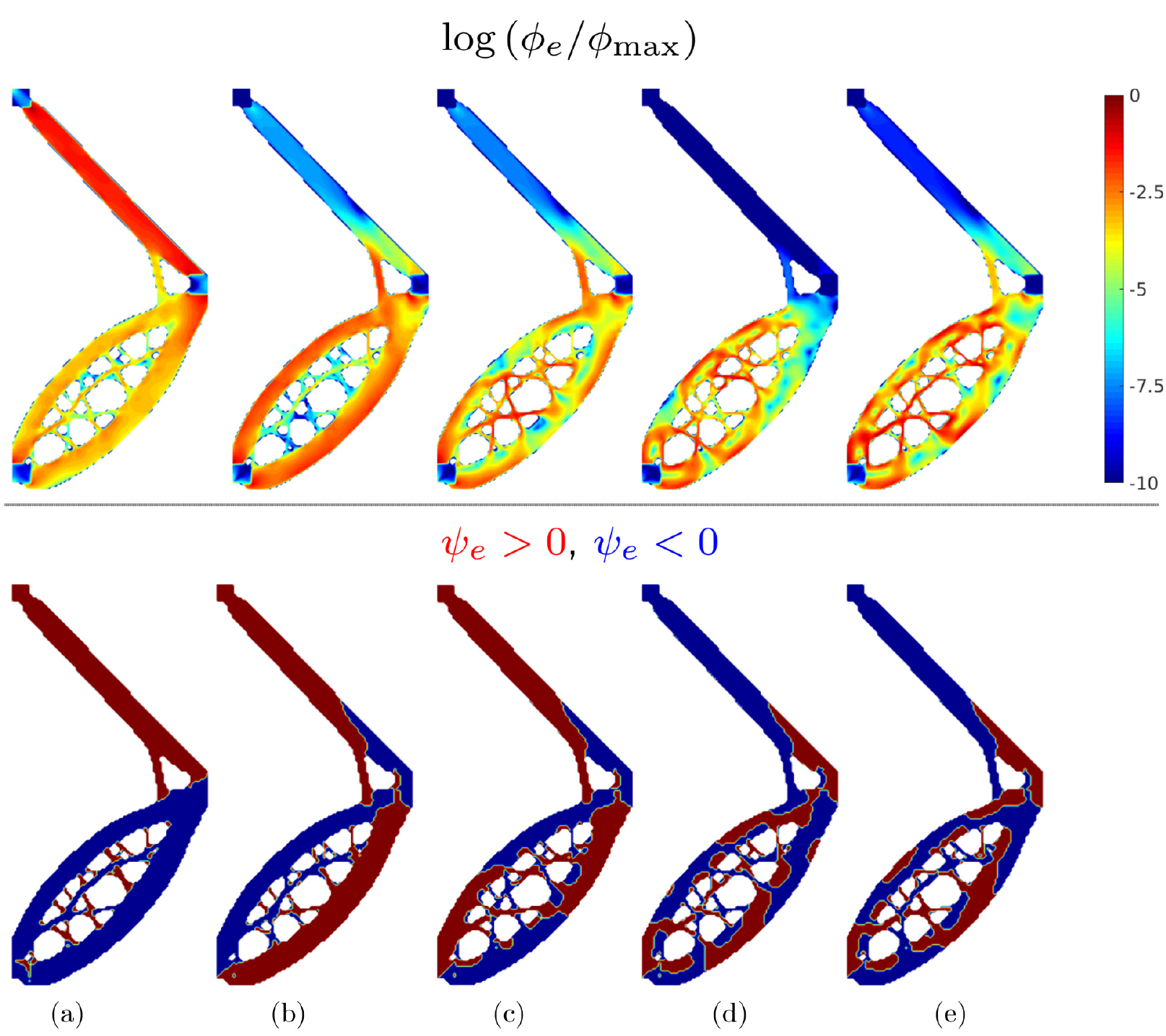}
 \caption{Logarithmic plot of the relative strain energy density (first row) and qualitative plot of the regions subjected to stiffening and softening (second row) for the design of \autoref{fig:OptimizedDesignsIncreasingPc} (g). (a) refers to the elastic displacement and (b--e) to the lowest four buckling modes}
 \label{fig:ComparisonBucklingModes}
\end{figure}

At each step we solve \eqref{eq:LinearizedEquilibriumEquation}, and compute the first 24 eigenpairs from \eqref{eq:GEPforMu}. A criterion based on the strain energy density, as described by \citet{gao-ma_15a}, is used for assessing the physical meaningfulness of the modes, i.e. for filtering out the ones corresponding to artificial deformations of low density regions. However, we point out that in the following testing we never experienced this issue. The lowest 12 buckling load factors are then included as constraints, in the fashion of \eqref{eq:EigConstraintMU}. For all the following examples we set $\alpha = 1.01$, for imposing the eigenvalue gaps.

The design variables are updated by the Method of Moving Asymptotes (MMA) \citep{svanberg_87a} and the optimization is run for a fixed number of 700 iterations.

The compliance and the fundamental buckling load for the design with uniform material distribution are $J_{0} = 4.714 \cdot 10^{-3}$ and $\lambda_{1, 0} = 0.662$, respectively.

The compliance design corresponding to a volume fraction of $f = 0.2$ and without any account for buckling consists of two identical bars, the upper one mainly in tension and the lower one in compression (see \autoref{fig:GeometricalSettingAndComplianceDesing} (b)). The strain energy density (defined on the element level as $\phi_{e} = \boldsymbol{\varphi}^{T}_{ie}\mathbf{k}_{e}\boldsymbol{\varphi}_{ie}$) corresponding to the lowest four buckling displacements is shown in \autoref{fig:GeometricalSettingAndComplianceDesing} (c) and we recognize global bending modes in the lower bar with progressively shorter wave lengths.

\begin{table*}[t]
 \caption{Compliance reduction factors and eigenvalue separation factors (as defined by \eqref{eq:DeltaDefinition}) for the optimized designs corresponding to increasing values of $\overline{P_{c}}$}
 \label{tab:ResultsComplianceAndBucklingLoad}
 \centering
  \begin{tabular}{l|c|ccccc}
   \hline\noalign{\smallskip}
   $\overline{P_{c}}$ & $J_{n}$ & \multicolumn{5}{c}{$\delta_{i}$} \\
   \noalign{\smallskip}\hline
   0.0  & 0.5564 & 1.69 & 3.48 & 5.16 & 6.57 & 7.63 \\
   0.5  & 0.5727 & 0.53 & 1.89 & 2.55 & 3.18 & 3.29 \\
   0.75 & 0.5903 & 6.1$\cdot 10^{-4}$ & 0.65 & 1.07 & 1.42 & 1.74 \\
   1.0  & 0.6241 & 6.6$\cdot 10^{-5}$ & 0.25 & 0.75 & 0.84 & 0.93 \\
   1.25 & 0.6899 & 5.4$\cdot 10^{-5}$ & 1.7$\cdot 10^{-4}$ & 2.6$\cdot 10^{-2}$ & 2.8$\cdot 10^{-2}$ & 8.3$\cdot 10^{-2}$ \\
   1.5  & 0.7556 & 6.9$\cdot 10^{-6}$ & 8.5$\cdot 10^{-5}$ & 1.4$\cdot 10^{-4}$ & 1.3$\cdot 10^{-1}$ & 1.4$\cdot 10^{-1}$ \\
   1.75 & 0.8311 & 3.7$\cdot 10^{-6}$ & 2.0$\cdot 10^{-5}$ & 3.7$\cdot 10^{-4}$ & 9.0$\cdot 10^{-4}$ & 2.3$\cdot 10^{-3}$ \\
   1.8  & 0.8664 & 3.1$\cdot 10^{-6}$ & 2.9$\cdot 10^{-5}$ & 3.0$\cdot 10^{-4}$ & 3.6$\cdot 10^{-4}$ & 4.4$\cdot 10^{-4}$ \\
   2.0  & 0.9729 & 1.1$\cdot 10^{-6}$ & 1.2$\cdot 10^{-6}$ & 1.6$\cdot 10^{-6}$ & 3.8$\cdot 10^{-5}$ & 5.6$\cdot 10^{-4}$ \\
   2.0$^{\ast}$  & 0.9648 & 1.5$\cdot 10^{-6}$ & 3.6$\cdot 10^{-6}$ & 5.8$\cdot 10^{-6}$ & 1$\cdot 10^{-5}$ & 1.5$\cdot 10^{-5}$ \\
   2.0$^{\ast\ast}$  & 0.9219 & 1.3$\cdot 10^{-6}$ & 5.9$\cdot 10^{-6}$ & 3.2$\cdot 10^{-5}$ & 6$\cdot 10^{-5}$ & 8.5$\cdot 10^{-5}$ \\
   \noalign{\smallskip}\hline
  \end{tabular}
\end{table*}

\subsection{Influence of the buckling constraint}
 \label{sSec:InfluenceBucklingConstraint}

The designs obtained for increasing values of $\overline{P_{c}}$ are shown in \autoref{fig:OptimizedDesignsIncreasingPc}. As the buckling constraint value is increased the upper tension bar becomes thinner and the material is relocated in the lower part of the domain, where a progressively wider truss--like configuration develops. This improves the bending stiffness, and therefore the buckling resistance, of the lower bar.

We also notice the appearance of a bar connecting the original upper and lower bars near point $c$. This element, becoming progressively thicker, has an important stiffening effect, ensuring a global behavior of the structure against the rotation around point $c$, which is highly strained in the fundamental buckling mode (see \autoref{fig:GeometricalSettingAndComplianceDesing} (c)).

The qualitative behavior of an optimized structure, e.g. the one corresponding to $\overline{P_{c}} = 1.8$, can be studied by looking at \autoref{fig:ComparisonBucklingModes}. Here the plots on the first row show the distribution of the strain energy densities for the elastic displacement and for the lowest four buckling modes. The plots in the second row display a subdivision of the regions where the stress energy density associated with a given displacement vector $\boldsymbol{v}$ and defined on the element level as $\psi_{e} = \boldsymbol{v}^{T}_{e}\mathbf{g}_{0}\boldsymbol{v}_{e}$, is either positive or negative. Positive (red) indicates stiffened domains and negative (blue) indicates softened domains.

For the pre--buckling displacement the deformation energy is highest in the upper bar, which is in traction, while the lowest part, mainly in compression, is softened. However, some of the inner bars are actually in tension, as they counteract the spreading of the two outer elements. As the load multiplier increases, the stiffness of the lower part is reduced, up to the occurrence of buckling. We can see from the strain energy distribution, and from the distribution of traction and compression regions, that the fundamental buckling mode consist of a global bending of the lower part. Moreover, the short connection bar near point $c$ is highly strained, which confirms its usefulness in restraining this instability mode, whereas the thin bars of the infill are only marginally affected by this deformation. For the higher buckling modes the deformation progressively affects also these and the modes become more localized.

\begin{figure}[t]
 \centering
  \includegraphics[scale = 0.5, keepaspectratio]
   {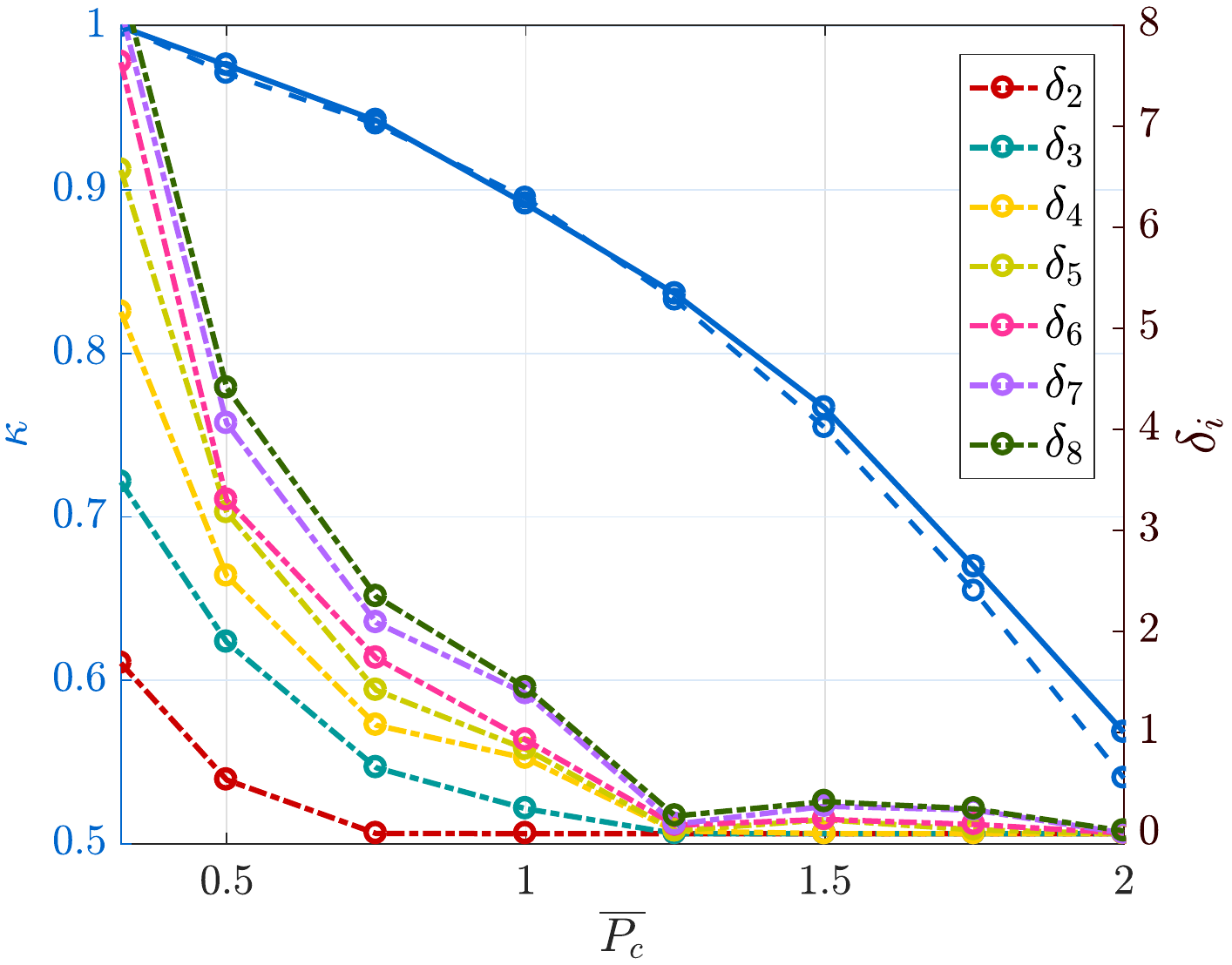}
  \caption{Performance of the designs optimized for increasing values of the lower bound on the fundamental buckling load factor. The blue curves show the trend of the stiffness ratio $\kappa$ when using $\mathcal{Q}_{4}$ (solid curve) and $\mathcal{Q}_{6}$ (dashed curve) elements for the discretization}
 \label{fig:ParetoSetUniform}
\end{figure}

Interestingly, some bending waves similar to those appearing for the compliance design can still be qualitatively recognized, especially if we look at the distribution of the stress energy $\psi_{e}$.

For the high value of the buckling constraint $\overline{P_{c}} = 2.0$ another linking bar appears (see \autoref{fig:OptimizedDesignsIncreasingPc} (h)). This element, however, turns out to configure an inconvenient deployment of material and this design clearly constitute a sub--optimal solution.

The introduction of the buckling constraint considerably increases the complexity of the feasible set and the chances of ending in a local minimum are increased as this constraint becomes more influential. A natural precaution can be to start with a non--uniform material distribution \citep{book:bendsoe-sigmund_2004}, closer to the design looked for.

The designs in \autoref{fig:OptimizedDesignsIncreasingPc} (i, j), for example, are obtained for the buckling constraint $\overline{P_{c}} = 2.0$ starting from the designs corresponding to $\overline{P_{c}} = 0$ and $\overline{P_{c}} = 1.0$, respectively. We see that the longer linking element is first made shorter and then removed, for the benefit of a slightly thicker upper bar. It is also interesting to notice that although the short bars in the infill have different configurations for the three last designs, this seems not to affect significantly the buckling response, since the critical load and the corresponding mode remain practically unchanged.

The performance of the obtained designs is quantitatively assessed by looking at the ratio between the final and initial compliance $J_{n} = J_{\rm f} / J_{0}$ and at the parameters
\begin{equation}
 \label{eq:DeltaDefinition}
 \delta_{i} = \lambda_{i}/\lambda_{1} - \alpha^{(i-1)} \: , \ i > 2
\end{equation}
which approach zero as the $i$--th buckling mode tends to coalesce on the fundamental one (with a separation of $\alpha^{(i-1)}$).

Numerical values of these parameters are collected in \autoref{tab:ResultsComplianceAndBucklingLoad} and their ``at--a--glance'' representation is given by \autoref{fig:ParetoSetUniform}. Here, the blue continuous curve, plotted against the left axis, shows the ratio between the stiffness of a design obtained for a given buckling constraint and that for $\overline{P_{c}} = 0$, namely $\kappa = J_{n,\overline{P_{c}}}/J_{n,\overline{P_{c}}=0}$. From the trend of $\kappa$ we clearly see that the designs optimized for a fixed volume fraction cannot retain the same stiffness as $\overline{P_{c}}$ is increased. Moreover, $\kappa$ becomes steeper after $\overline{P_{c}} \approx 1.5$ and this indicates that for higher values of the buckling constraint a very compliant structure emerges during the optimization, and even convergence to a well--connected design can be hampered.

\begin{figure}[tb]
 \centering
  \includegraphics[scale = 0.7, keepaspectratio]
   {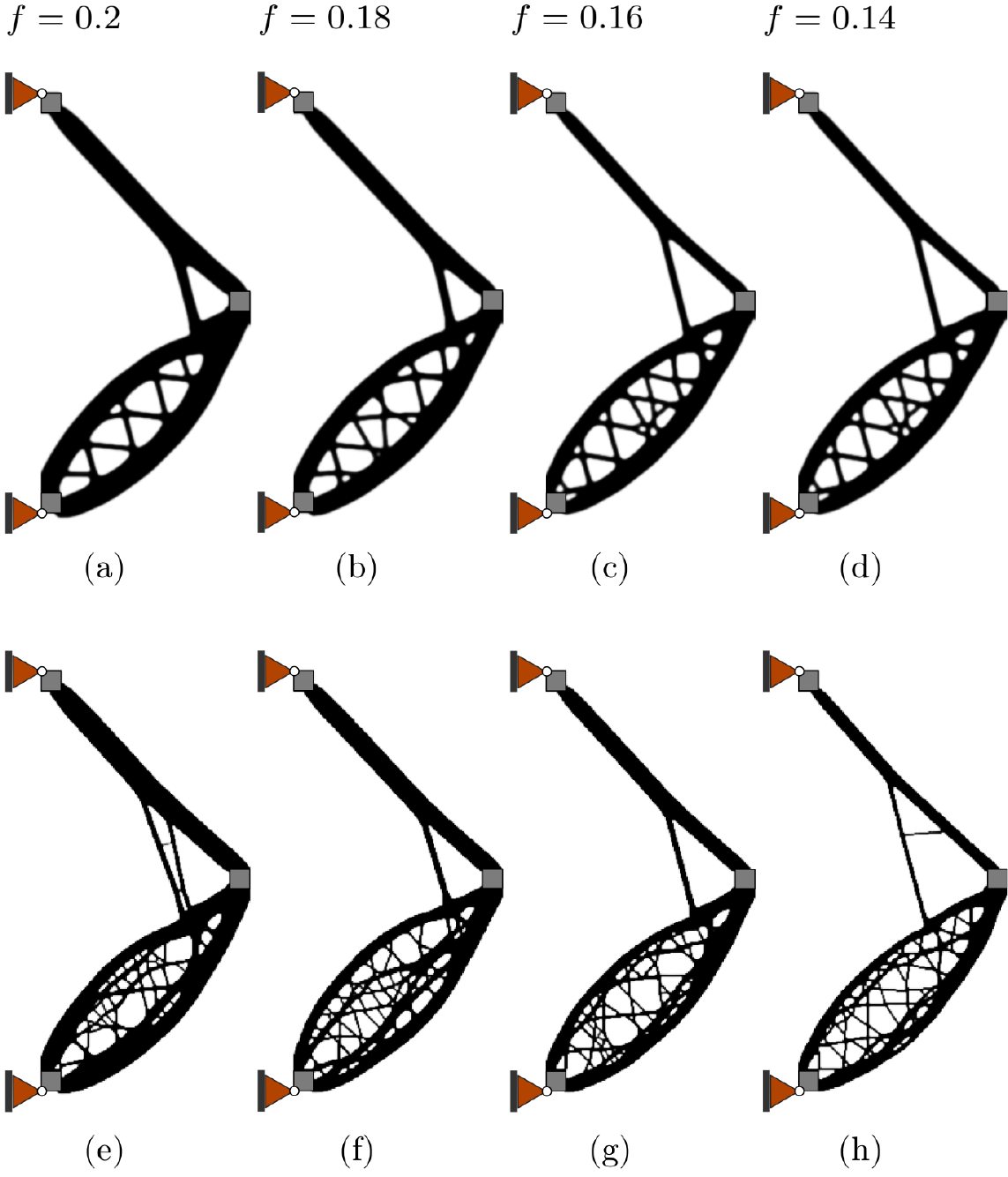}
 \caption{Designs obtained on the mesh ${\rm \Omega}_{2} = 180\times 420$ for decreasing values of the volume constraint and for the filter sizes of $r_{\rm min} = 4$ elements (first row) and $r_{\rm min} = 2$ elements (second row). The initial guess is, for all the cases, the corresponding compliance design}
 \label{fig:DesignsOnFinerMesh}
\end{figure}

The behavior of the eigenvalue separation parameters defined in \eqref{eq:DeltaDefinition} is shown by the dashed curves plotted against the right axis and we acknowledge the progressive activation of an increasing number of buckling constraints. Already for $\overline{P_{c}} = 0.75$ we have a bi--modal $\lambda_{1}$ and for $\overline{P_{c}} = 1.25$ the lowest eight eigenvalues are all very close. Then, they become all simultaneously active for the highest value $\overline{P_{c}} = 2.0$.

The general trend we recognize is that the finer the bar distribution in the infill the more the buckling modes group together. From this simple example we can see that a large number of buckling modes may contribute to the optimization, increasing the complexity of the problem, and this is bound to get even more pronounced for higher mesh resolutions and lower volume fractions.

\begin{figure*}[tb]
 \centering
  \subfloat[]{\includegraphics[scale = 0.475, keepaspectratio]
   {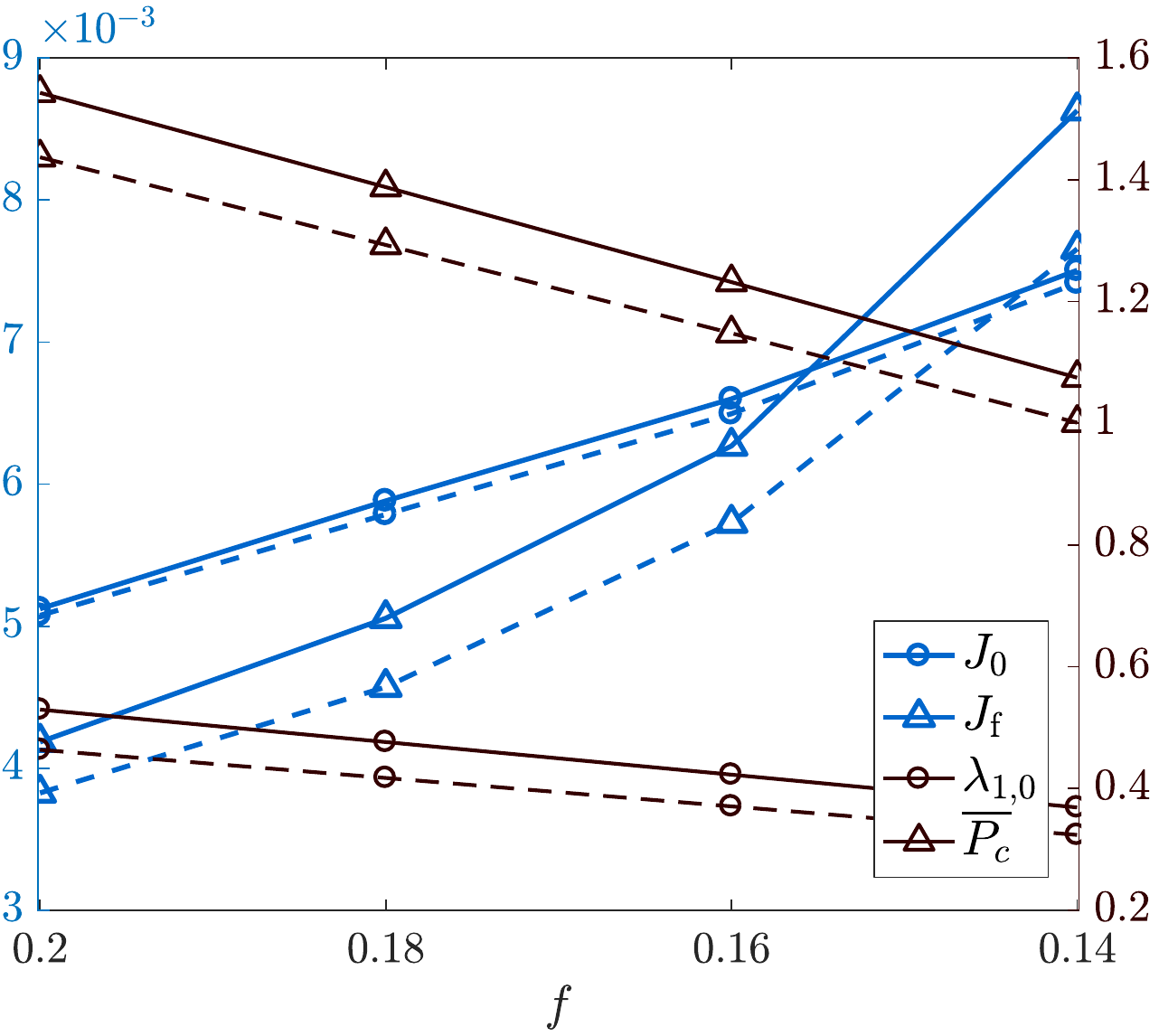}} \qquad
  \subfloat[]{\includegraphics[scale = 0.475, keepaspectratio]
   {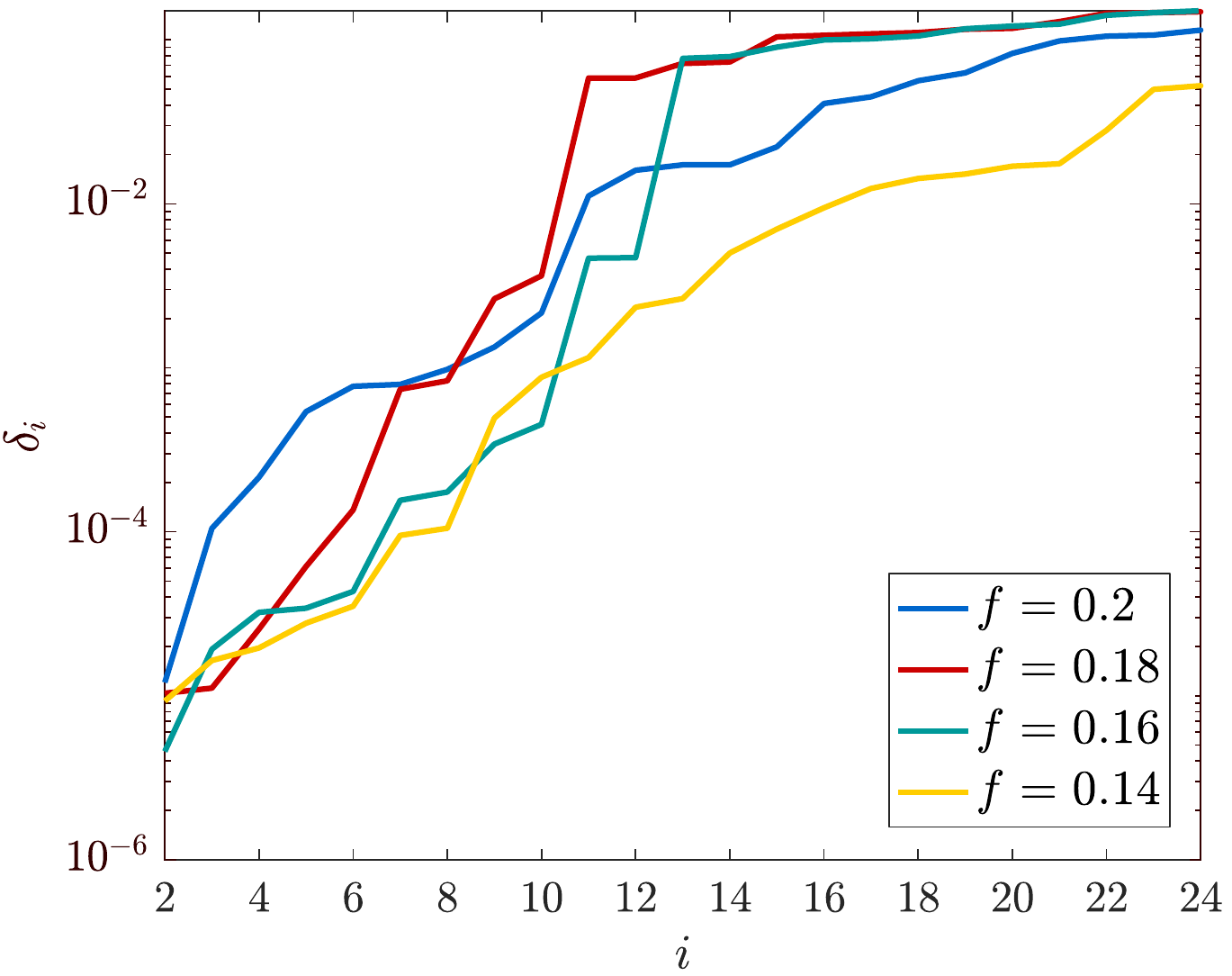}}
 \caption{(a) shows the modification of the compliance (referred to left axis) and of the fundamental buckling load factor (referred to the right axis) for the uniform and optimized design on the finer mesh ${\rm \Omega}_{2}$. Solid lines refer to the case $r_{\rm min} = 4$ and dashed ones to the case $r_{\rm min} = 2$. (b) shows the trend of the eigenvalue gap parameters $\delta_{i}$ defined in \eqref{eq:DeltaDefinition} for the designs with low filter size and decreasing allowed volume}
 \label{fig:FinerMeshBehavior}
\end{figure*}

\subsection{Designs obtained on a finer discretization}
 \label{sSec:FinerMeshDesigns}

\autoref{fig:DesignsOnFinerMesh} displays some designs obtained on a finer discretization ${\rm \Omega}_{2} = 180 \times 420$, for decreasing values of the volume constraint.

The designs in the first row refer to a filter size scaled according to the mesh size (i.e. $r_{\min} = 4$ elements), while those in the second row are obtained for the same filter size (compared to element size) used for the coarser mesh ${\rm \Omega}_{1}$ (i. e. $r_{\min} = 2$ elements). In all the cases, the corresponding compliance design is used as initial guess, to facilitate the obtainment of a good solution.

Now 32 eigenpairs are computed and the lowest 24 buckling loads are included as constraints, since we expect a higher number of modes to coalesce.

All these designs ideally relate to the one corresponding to $\overline{P_{c}} = 1.8$, shown in \autoref{fig:OptimizedDesignsIncreasingPc} (g). However, since the approximate buckling load reduces, both due to the mesh refinement and as the volume is reduced, it becomes extremely hard for the optimizer to find a design coping with such a high buckling constraint. Therefore $\overline{P_{c}}$ is scaled according to the trend shown in \autoref{fig:FinerMeshBehavior} (a), and its numerical values are $1.55$, $1.3$, $1.148$ and $1.002$, for the four decreasing volume fractions (found by numerical experiments). From this figure we see that the buckling load factor ($\lambda_{1,0}$) and compliance ($J_{0}$) of the uniform material distribution are linearly reduced and increased, respectively, as $f$ decreases. On the other hand, the compliance of the optimized design ($J_{\rm f}$) seems to increase more than proportionally. The same trend is recognized for both the cases of filter radii.

\begin{figure*}[t]
 \centering
  \includegraphics[scale = 0.7, keepaspectratio]
   {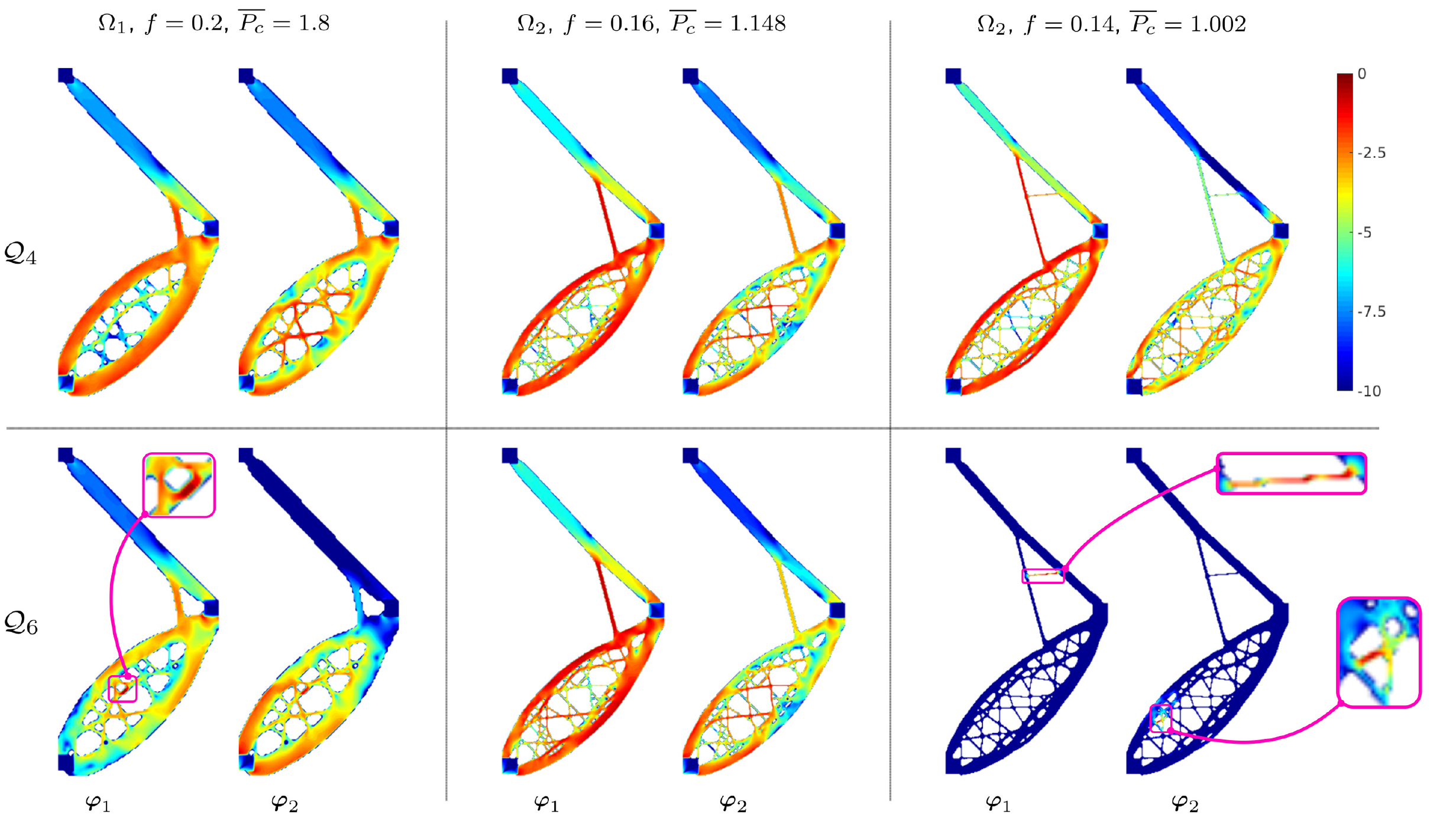}
 \caption{Comparison of the strain energy distribution ($\log(\phi_{e}/\phi_{\rm max})$) of the two lowest buckling modes for three optimized designs, analyzed with $\mathcal{Q}_{4}$ and $\mathcal{Q}_{6}$ elements. We observe that in all the cases the buckling mode predicted by $\mathcal{Q}_{6}$ elements is more localized, involving some fine bar. As an extreme case, for the design obtained on the finer mesh ${\rm \Omega}_{2}$ and for $ f  = 0.14$ and $r_{\rm min} = 2$, the fundamental buckling mode involves an isolated member}
 \label{fig:Q4Q6ModesDifference}
\end{figure*}

\begin{figure*}[tb]
 \centering
  \subfloat[$\overline{P_{c}} = 1.5$]{
   \includegraphics[scale = 0.215, keepaspectratio]
   {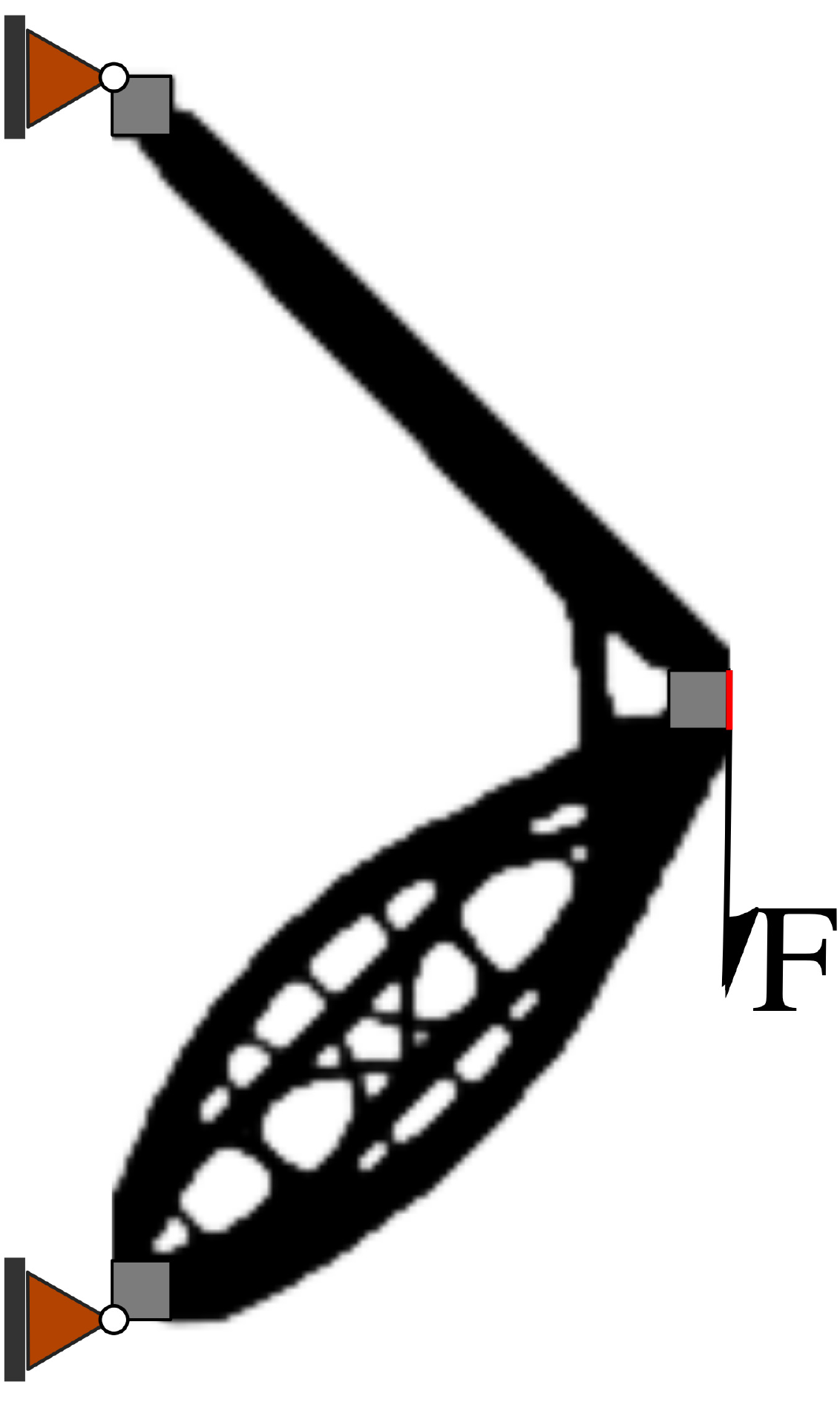}} \qquad
  \subfloat[$\overline{P_{c}} = 1.8$]{
   \includegraphics[scale = 0.215, keepaspectratio]
   {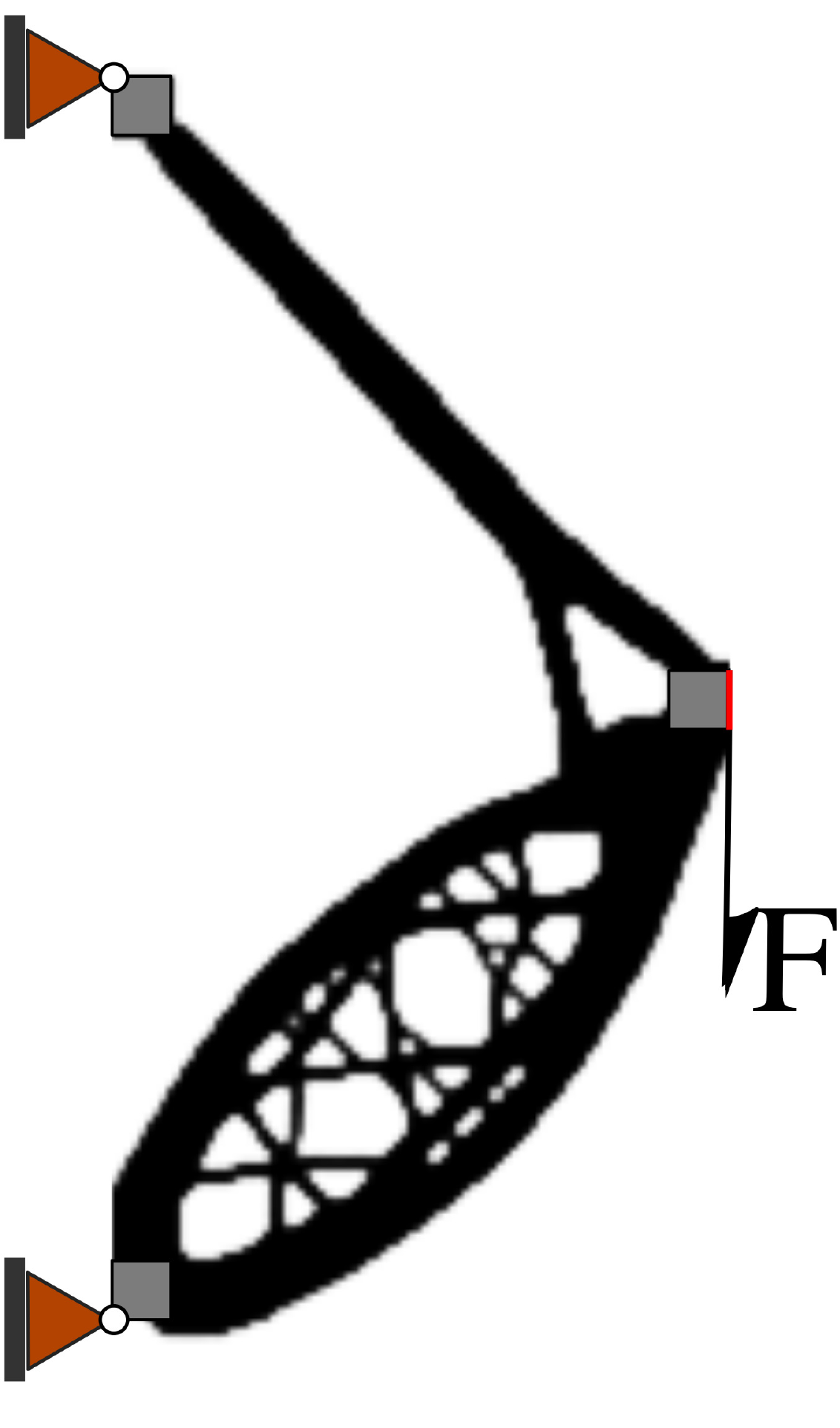}} \qquad
  \subfloat[$\overline{P_{c}} = 2.0$]{
   \includegraphics[scale = 0.215, keepaspectratio]
   {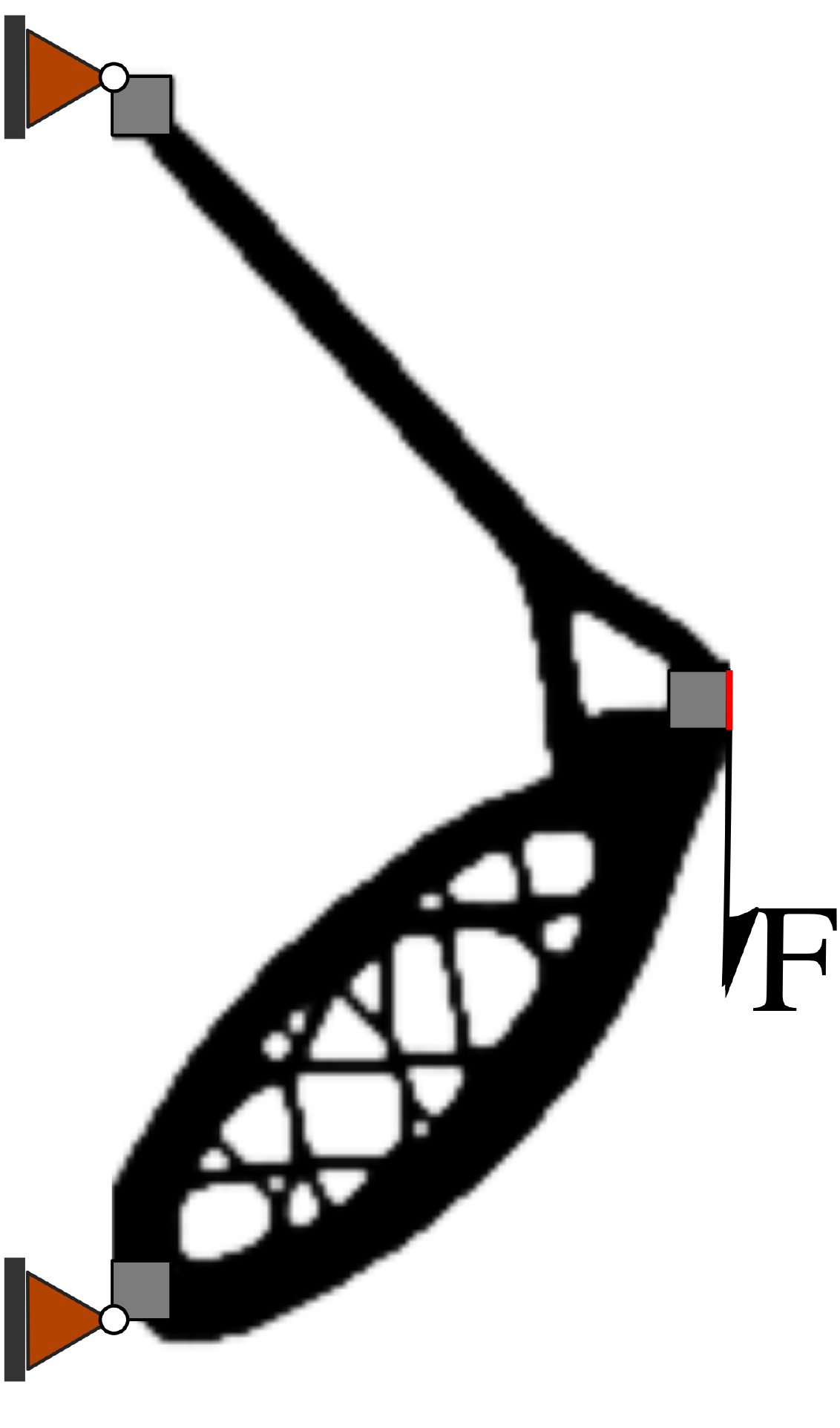}} \qquad
  \subfloat[${\rm\Omega}_{2}$,$ f  = 0.14$]{
   \includegraphics[scale = 0.215, keepaspectratio]
   {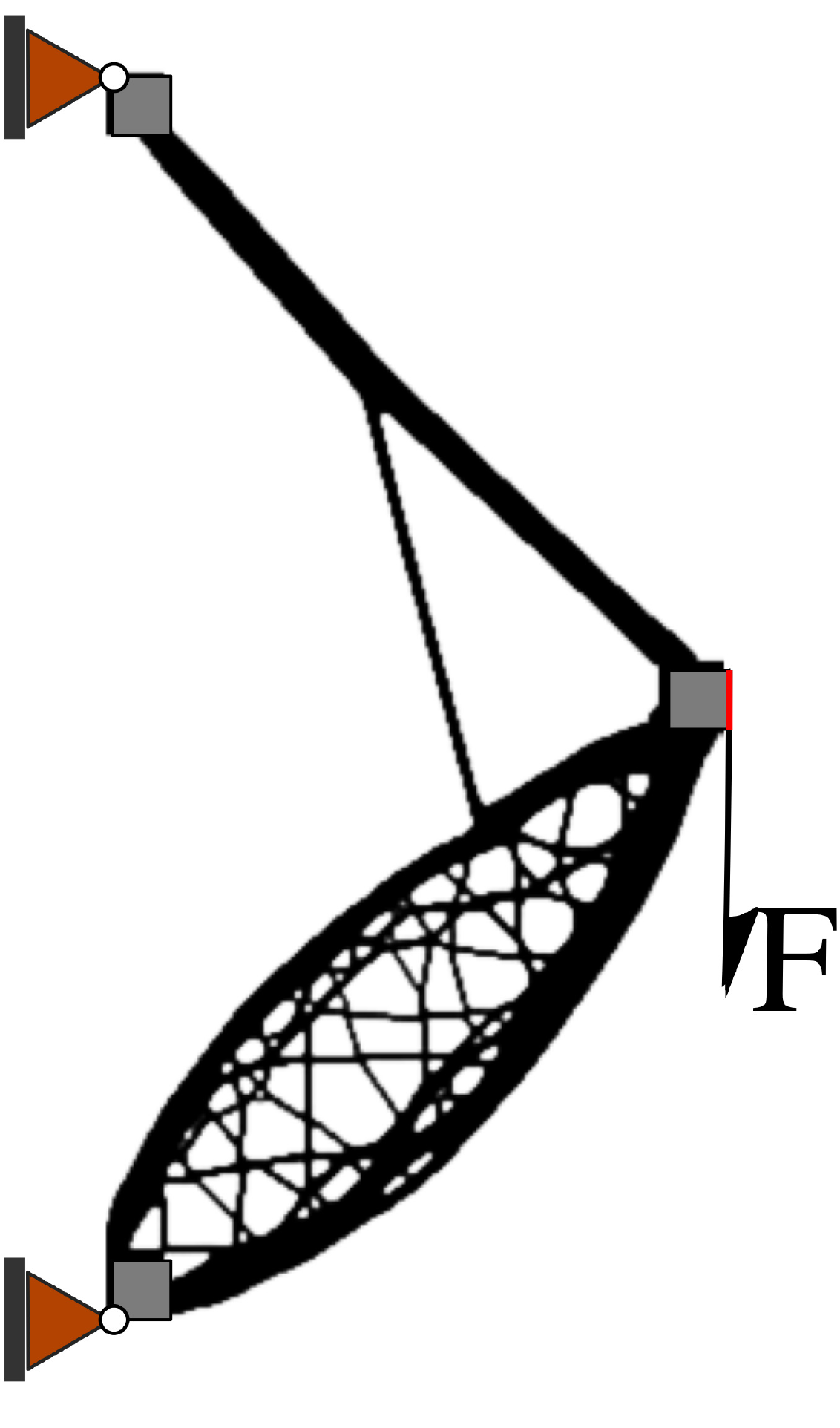}} \qquad
  \subfloat[]{
   \includegraphics[scale = 0.36, keepaspectratio]
   {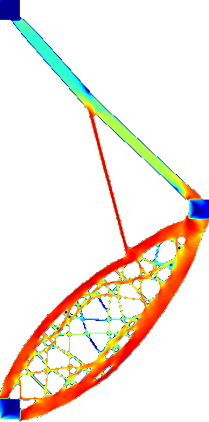}}
 \caption{Designs obtained by using $\mathcal{Q}_{6}$ elements for the state discretization. (a--c) refers to the mesh ${\rm \Omega}_{1}$, and compares with designs of \autoref{fig:OptimizedDesignsIncreasingPc}. (d) compares with the design of \autoref{fig:DesignsOnFinerMesh} (h) and (e) shows the distribution of the strain energy density corresponding to its fundamental buckling mode}
 \label{fig:OptimizedDesignsQ6}
\end{figure*}

As the volume fraction $f$ is reduced, a more complicated infill develops in the lower part, with the bars becoming thinner. The connection bar near the right end also becomes thinner and moves progressively away from point $c$. These modifications are further emphasized for designs obtained with the smaller filter size, where finer structural features are allowed.

\autoref{fig:FinerMeshBehavior} (b), shows the trend of the $\delta_{i}$ parameters, as defined by \eqref{eq:DeltaDefinition}, corresponding to the constrained eigenvalues. We see that more modes are becoming active, as expected. For the volume fraction $f  = 0.16$ the first ten $\delta_{i}$ are below $10^{-2}$, which means that $\lambda_{10}$ is only $12\%$ larger than $\lambda_{1}$. For the case $f  = 0.14$ the same ten go even below $10^{-3}$, and all the $\delta_{i}$ are below $5\cdot 10^{-2}$, which means that $\lambda_{24} < 1.3 \lambda_{1}$.

\begin{figure}[tb]
 \centering
  \subfloat[$\overline{P_{c}} = 0.5$]{
   \includegraphics[scale = 0.2, keepaspectratio]
   {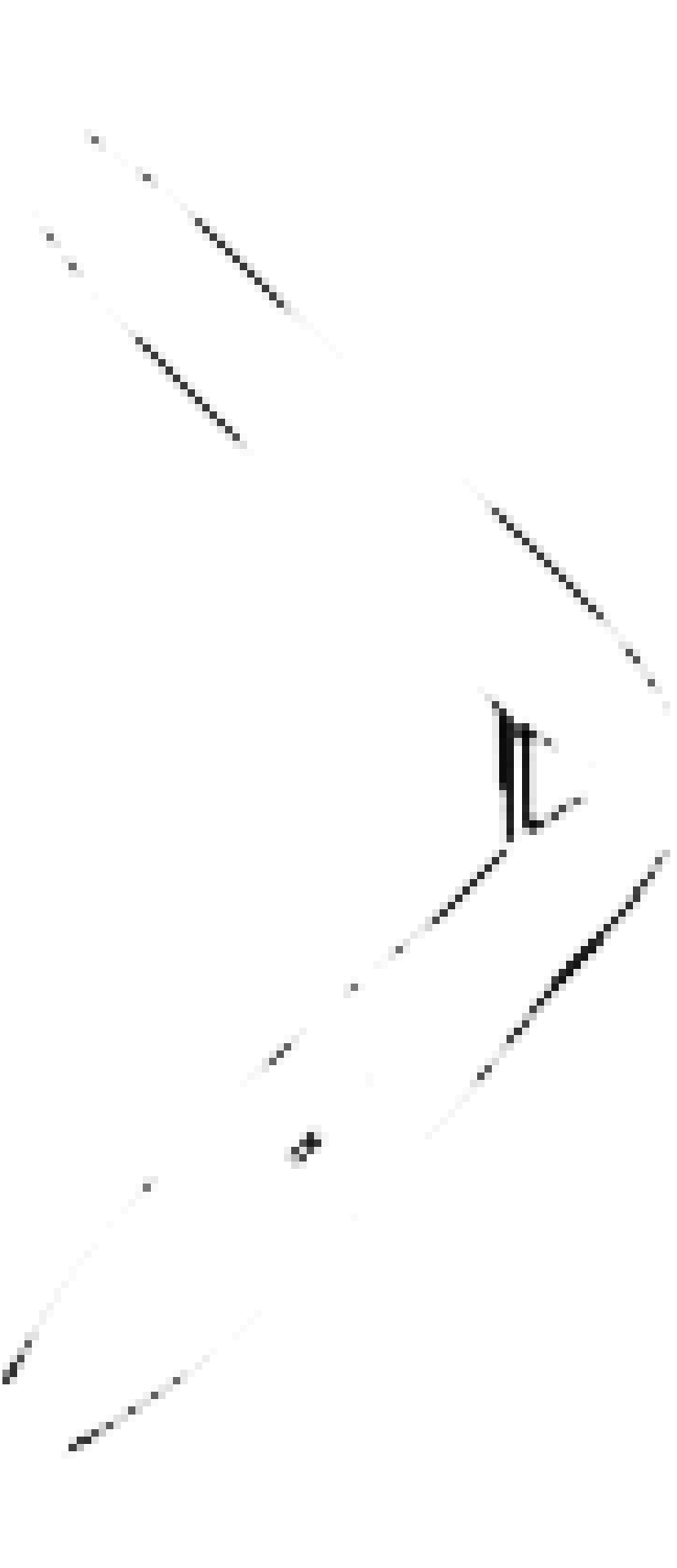}} \quad
  \subfloat[$\overline{P_{c}} = 0.75$]{
   \includegraphics[scale = 0.2, keepaspectratio]
   {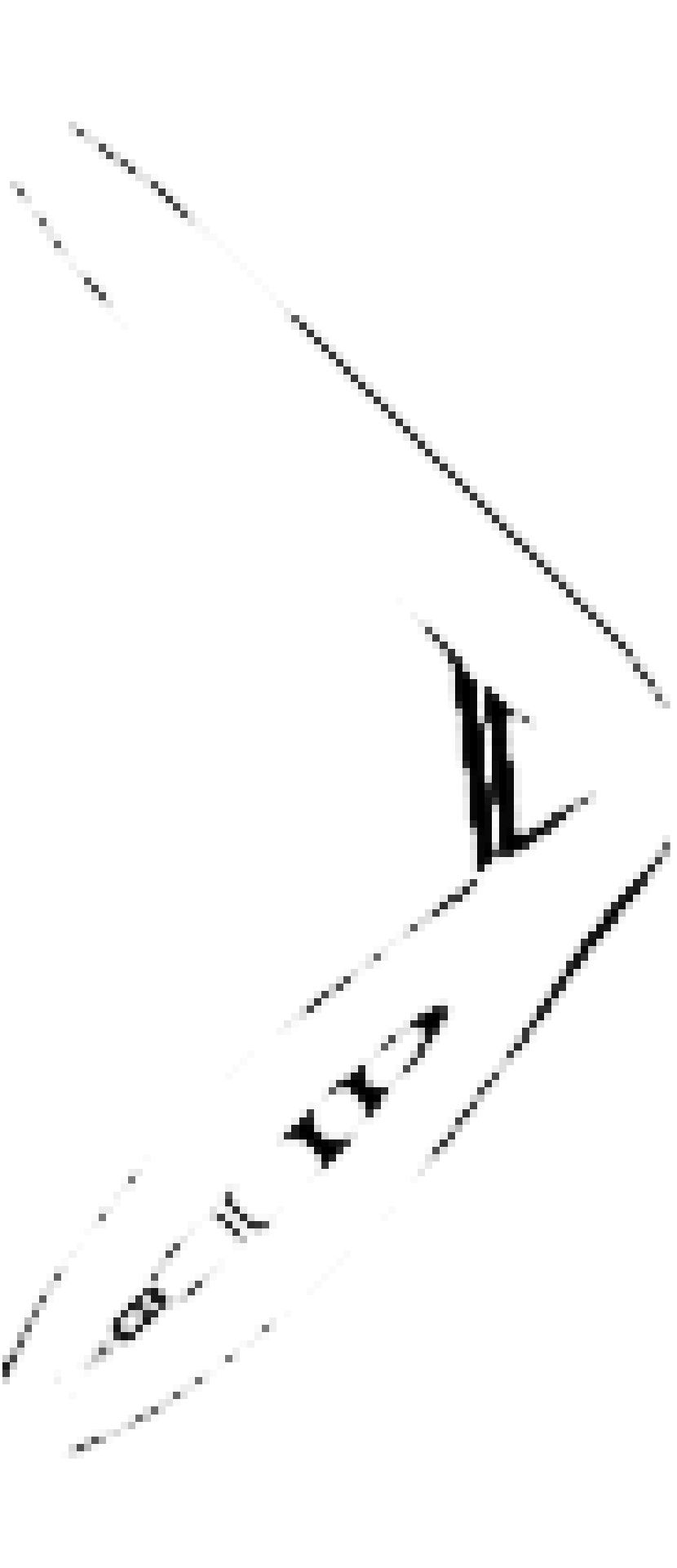}} \quad
  \subfloat[$\overline{P_{c}} = 1.0$]{
   \includegraphics[scale = 0.2, keepaspectratio]
   {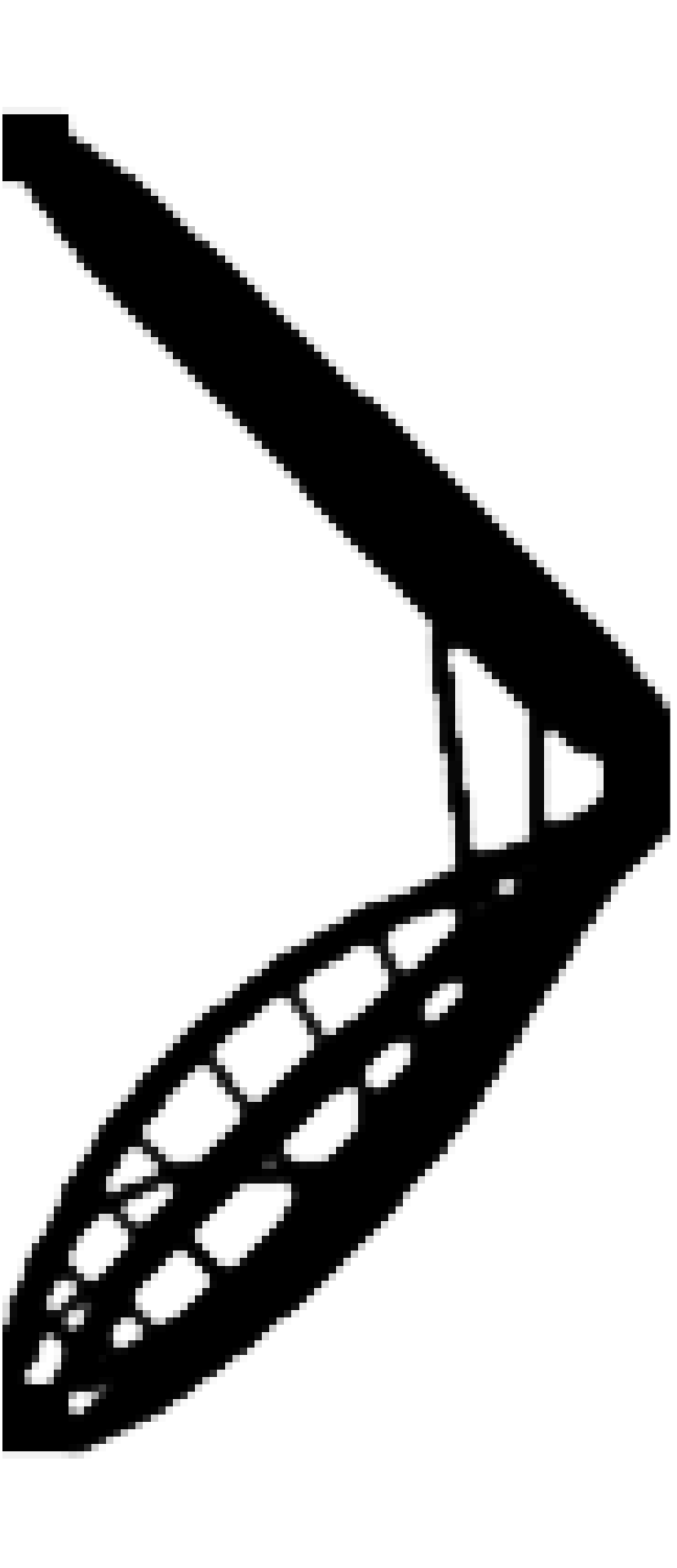}} \quad
  \subfloat[$\overline{P_{c}} = 2.0$]{
   \includegraphics[scale = 0.2, keepaspectratio]
   {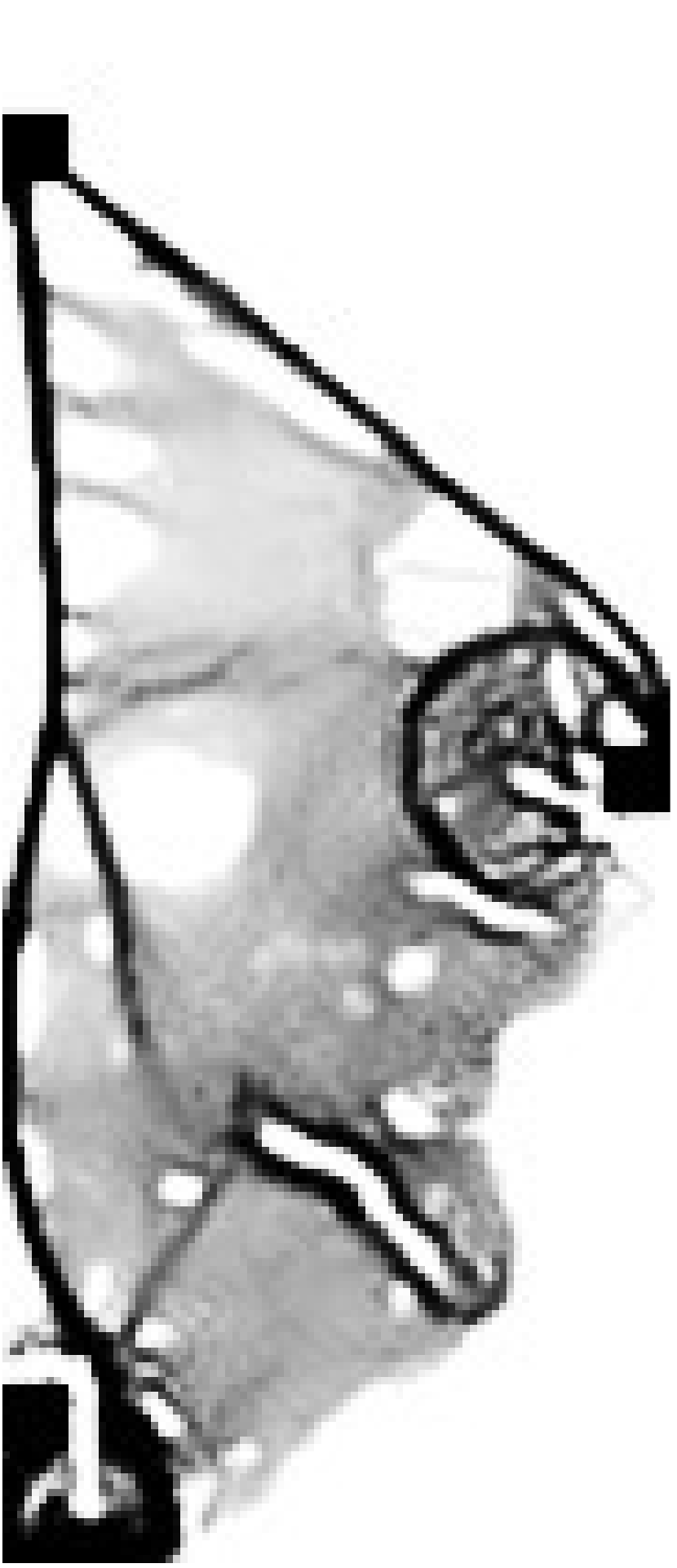}}
 \caption{Designs obtained by using an inconsistent sensitivity. (a) and (b) show the difference with respect to the corresponding designs of \autoref{fig:OptimizedDesignsIncreasingPc}. For these we have the compliance ratios $J_{n} = 0.5753$ and $J_{n} = 0.5963$. For the design of (c) we have $J_{n} = 0.6597$ and for the material distribution in (d) the optimization breaks at iteration 150}
 \label{fig:ComparisonDesingDifferences}
\end{figure}

\subsection{Analysis and design with $\mathcal{Q}_{6}$ elements}
 \label{sSec:ComparisonResultsQ6}

All the previous designs have been obtained by using $\mathcal{Q}_{4}$ elements for the discretization. We now analyze their response also by using $\mathcal{Q}_{6}$ elements.

From this analysis we find that the computed compliance and fundamental buckling load factor are very similar, with differences of the order of few \textperthousand\hphantom{} and below $1\%$, respectively. Nevertheless, important differences in the nature of the buckling modes can be seen by looking at \autoref{fig:Q4Q6ModesDifference}, where the strain energy density associated with the two lowest buckling modes, as predicted by $\mathcal{Q}_{4}$ and $\mathcal{Q}_{6}$ discretizations, is displayed for three representative designs.

For the design corresponding to $\overline{P_{c}} = 1.8$ and the mesh ${\rm \Omega}_{1}$, the buckling load predicted by the $\mathcal{Q}_{6}$ discretization is $\lambda_{1} = 1.791$, i.e. only $0.5\%$ lower than the one given by $\mathcal{Q}_{4}$ elements. However, the leftmost set of plots in \autoref{fig:Q4Q6ModesDifference} shows that the strain energy is much more localized in the infill, where one of the thin bars undergoes local buckling. For the second mode the maximum strain energy is still occurring in the same bar.

The central set of plots concerns the response of the design from \autoref{fig:DesignsOnFinerMesh} (g), and similar considerations can be made. Again, some of the bars in the infill appear to be highly strained for the buckling mode predicted by the $\mathcal{Q}_{6}$ elements, though the numerical value of the buckling load factor is only $0.6\%$ smaller.

As an extreme case, we include the rightmost set of plots in \autoref{fig:Q4Q6ModesDifference}, concerning the design from \autoref{fig:DesignsOnFinerMesh} (h). This one shows some very thin bars and we acknowledge the completely different nature of the buckling modes, which become extremely local as predicted by the $\mathcal{Q}_{6}$ elements. The fundamental mode involves only a thin bar, which clearly indicates a non--optimal feature, described by very few elements. Also the second buckling mode involves the failure of a small bar of the infill, which is subjected to a high compressive stress. These mechanisms are completely missed when using $\mathcal{Q}_{4}$ elements, where such regions seem to be not much strained.

At this point it should become apparent that even though the differences in the predicted buckling load factors appear to be negligible, the substantial discrepancy of the buckling modes is expected to definitely play a role in the sensitivity, and therefore in the optimization outcomes.

In order to investigate this, the previous examples have been solved again by using $\mathcal{Q}_{6}$ elements. Design variables are referred to the element center and therefore filtering is still needed for the stability of the approximation \citep{jog-haber_96a}.

Some of the resulting designs are shown in \autoref{fig:OptimizedDesignsQ6}, and we notice that the main differences with those collected in \autoref{fig:OptimizedDesignsIncreasingPc} concern, as expected, the bars of the infill, which now appear to be generally thicker.

Looking back at \autoref{fig:ParetoSetUniform}, we see the difference in the performance of these designs, compared to those obtained by using $\mathcal{Q}_{4}$ elements. As we already pointed out, the discretized structure is now softer and the buckling load of the uniform design is about $2\%$ lower than the one predicted by conforming $\mathcal{Q}_{4}$ elements. However, referring to the compliance ratio $J_{n}$ and to the stiffness parameter $\kappa$, the comparison is still consistent.

For low values of $\overline{P_{c}}$ the differences between the two sets of results are negligible. Then the curves start to drift as $\overline{P_{c}}$ is increased and the stiffness ratio $\kappa$ corresponding to the structures optimized with $\mathcal{Q}_{6}$ elements is systematically lower. This is again expected, as the softening of the discretization implies that the buckling constraint is becoming more demanding. The difference in the $\kappa$ value attains, for $\overline{P_{c}} = 2$, a maximum of $3\%$.

For the designs on the finer mesh ${\rm \Omega}_{2}$ we even recognize more evident changes of the topology, with the removal of weak thin bars. For example, \autoref{fig:OptimizedDesignsQ6} (d) shows the result corresponding to $ f  = 0.14$ and we notice the removal of the thin bar responsible for the buckling. From (e) we then see that the fundamental buckling mode now involves more globally the lower part. Again, the cross check of all these designs, i.e. their analysis by means of $\mathcal{Q}_{4}$ elements results in very little numerical differences for both the compliance and the buckling load.

In summary, from our experience, the use of $\mathcal{Q}_{6}$ in stead of $\mathcal{Q}_{4}$ elements seems to affect only little the value of the response parameters: compliance and buckling load factors. Also, only minor differences appear in optimized designs. However, the predicted buckling modes can be different from those given by $\mathcal{Q}_{4}$ discretizations, where the stiffness of small features may be overestimated. This is particularly pronounced for higher order buckling modes. Therefore, considering that the use of $\mathcal{Q}_{6}$ elements neither increases computational cost nor over--complicates implementation, $\mathcal{Q}_{6}$ elements are in general recommended but do not seem imperative for linear buckling topology optimization.

\begin{figure}[tb]
 \centering
  \includegraphics[scale = 0.675, keepaspectratio]
   {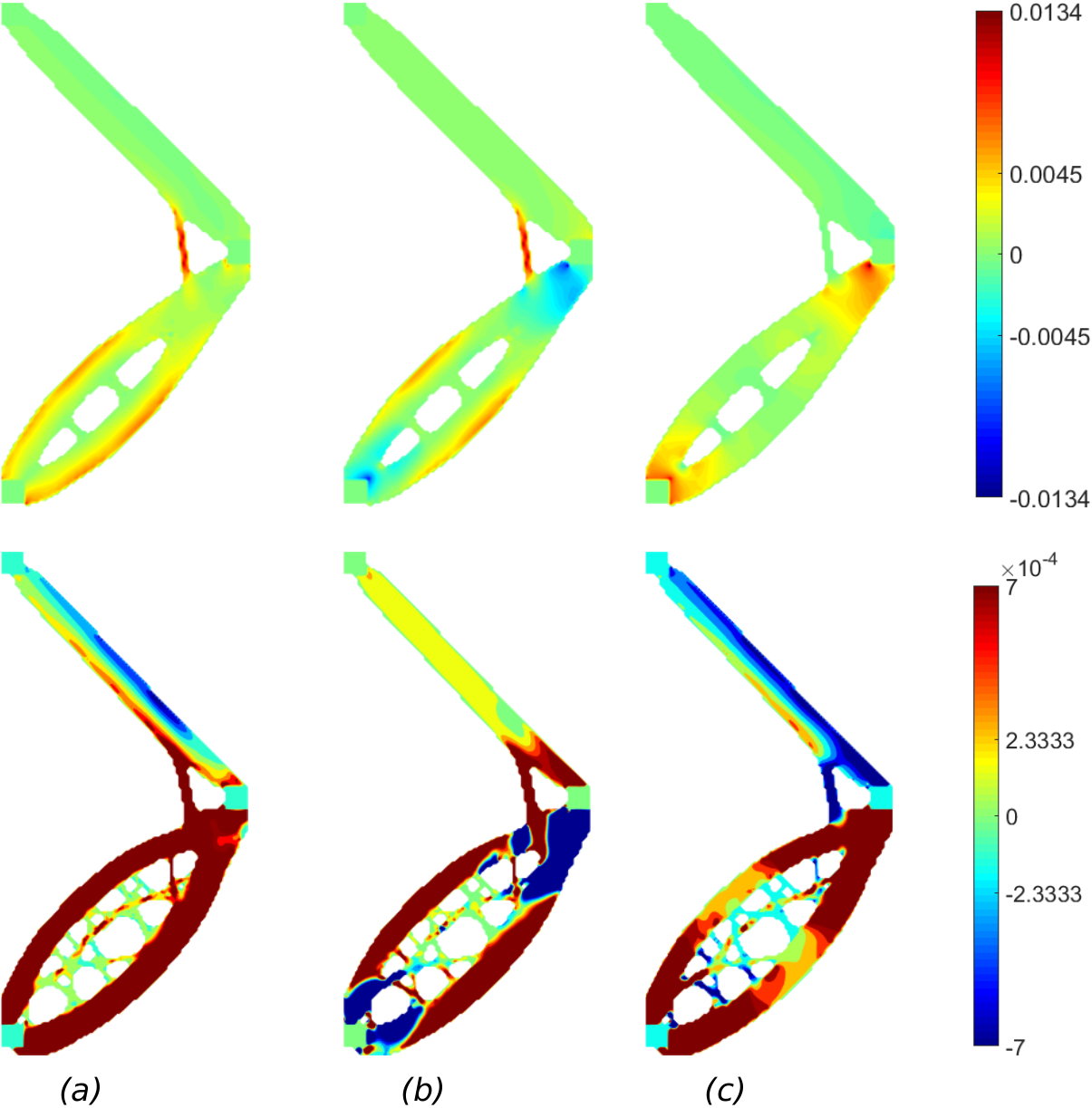}
 \caption{Plot of the sensitivity $\partial\lambda_{1}/\partial x_{e}$ (a) and of the frequency--like (b) and adjoint contributions (c) for the designs corresponding to $\overline{P_{c}} = 1.0$ (first row) and $\overline{P_{c}} = 1.8$ (second row). The color scale is referred to the same range, which is that of the total sensitivity}
 \label{fig:SensitivityPlotsDesignP1d0}
\end{figure}

\begin{figure*}[t]
 \centering
  \includegraphics[scale = 0.5, keepaspectratio]{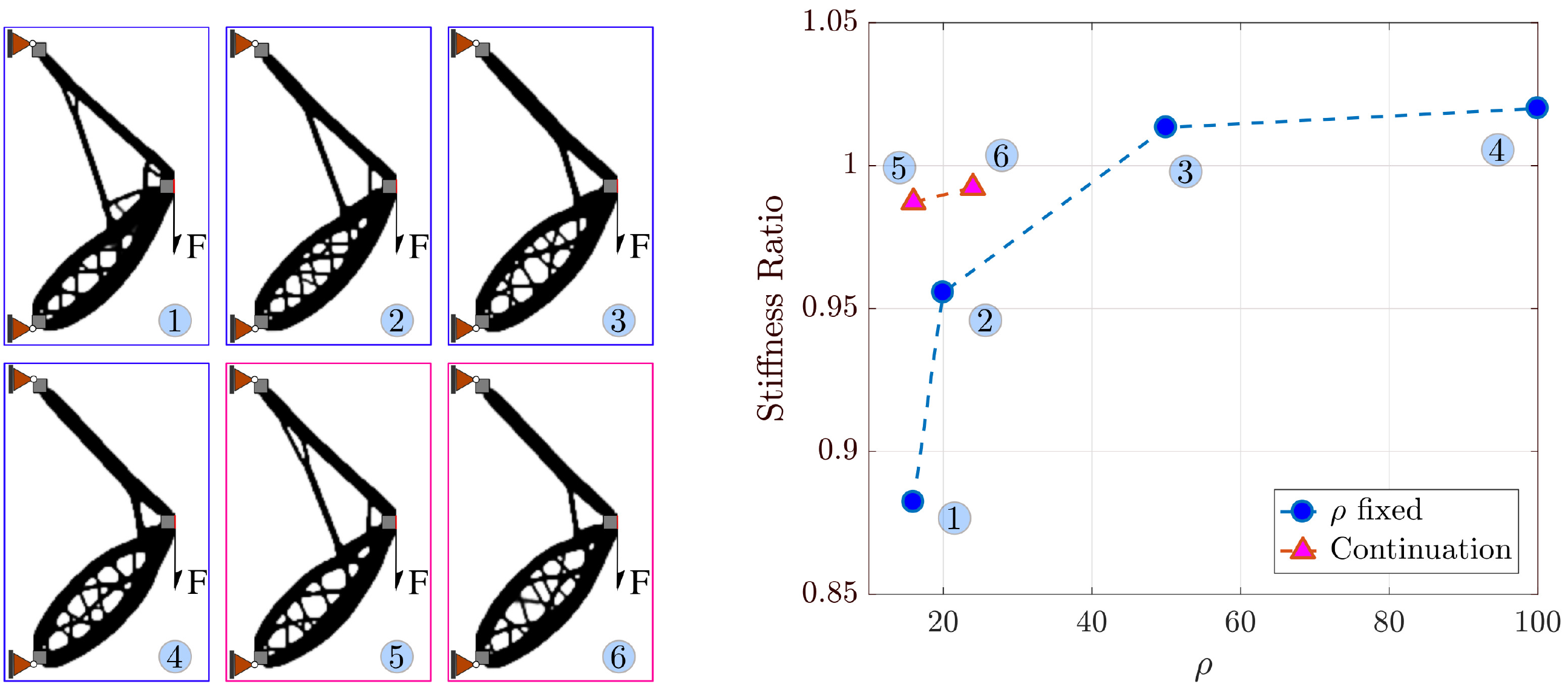}
 \caption{Designs obtained by using the KS function for aggregating the single constraints \eqref{eq:EigConstraintMU}, for the case $\overline{P_{c}} = 1.8$ and for several values of the aggregation parameter $\rho$. The graph on the right plots the ratio between the stiffness of these designs and that of the design in \autoref{fig:OptimizedDesignsIncreasingPc} (g)}
 \label{fig:KSaggregationComparison}
\end{figure*}

\subsection{Influence of the adjoint sensitivity term}
 \label{sSec:AdjointSensitivityTerm}

We now elaborate on the use of inconsistent sensitivities for the buckling constraint.

The set of examples discussed in \autoref{sSec:InfluenceBucklingConstraint} have been solved again considering only the frequency--like term of the sensitivity \eqref{eq:SensitivityLambda}, i.e.
\begin{equation}
 \left( \frac{\partial \lambda_{i}}{\partial x_{e}} \right)_{\rm inc.} = \boldsymbol{\varphi}^{T}_{i}\left( \frac{\partial \mathbf{K}}{\partial x_{e}} + \lambda_{i} \frac{\partial \mathbf{K}_{\sigma}}{\partial x_{e}} \right) \boldsymbol{\varphi}_{i}
 \ ,
\end{equation}
and neglecting the adjoint term $-\lambda_{i}\mathbf{v}^{T}\left[\partial \mathbf{K}/\partial x_{e}\right]\mathbf{u}_{0}$.

For low values of the buckling constraint ($\overline{P_{c}} = 0.5$ and $\overline{P_{c}} = 0.75$), this causes only marginal differences in the optimized design, mainly localized along the boundaries (see \autoref{fig:ComparisonDesingDifferences} (a, b)). However, we point out that already for these cases we have a (small) worsening of the performance, as $J_{n}$ increases by $0.4\%$ and $1\%$ respectively, compared to what reported in \autoref{tab:ResultsComplianceAndBucklingLoad}.

For the buckling constraint $\overline{P_{c}} = 1.0$ the differences become much more evident (see \autoref{fig:ComparisonDesingDifferences} (c)) and the increase in $J_{n}$ is about $6\%$. For higher values of the buckling constraint even the obtainment of a good design becomes difficult and, as an extreme situation, we can look at \autoref{fig:ComparisonDesingDifferences} (d), showing the material distribution at iteration 150 for the case $\overline{P_{c}} = 2.0$. Here the optimization clearly breaks down.

To further clarify the importance of the adjoint term in correctly driving the optimization, we refer to \autoref{fig:SensitivityPlotsDesignP1d0}, showing the distribution of the total sensitivity (a) and of the frequency--like (b) and adjoint (c) contributions for the optimized designs corresponding to $\overline{P_{c}} = 1.0$ and $\overline{P_{c}} = 1.8$ reported in \autoref{fig:OptimizedDesignsIncreasingPc}.

The distribution of the frequency--like term qualitatively resembles that of the total sensitivity, especially for the lower value of the buckling constraint. However, we notice some strongly negative values in the compressed regions, near zones where the load and displacements are prescribed. Therefore, this contribution alone would remove material from these regions (and this can be clearly recognized comparing the design of \autoref{fig:ComparisonDesingDifferences} (c) with the one in \autoref{fig:OptimizedDesignsIncreasingPc} (c)). This would lead to a very unstable design for high values of $\overline{P_{c}}$, as this tendency becomes more prominent.

On the other hand, the adjoint term attains its maximum values precisely in those regions, which are subjected to high stress gradients in the pre--buckling state. We also see that the magnitude of the adjoint contribution is definitely comparable to the frequency--like one and therefore their sum restore a consistent sensitivity distribution, preventing the material to be removed from highly stressed regions (see \autoref{fig:SensitivityPlotsDesignP1d0} (a)).

The discussion above should clearly show that the use of an inconsistent sensitivity for topology optimization problems involving linearized buckling, sometimes appearing in the literature \citep{munk-etal_16a}, should be avoided. In fact, this is at least the source of sub--optimal results and, for high values of the buckling constraints, may seriously hamper the whole optimization process.

\begin{table}[tb]
 \caption{Results obtained with the KS function for several values of the aggregation parameter. The last two rows refer to the application of a continuation scheme on the aggregation parameter, e.g. starting from $\rho=16$ and doubling it each 100 iterations up to $\rho = 128$}
 \label{tab:ResultsKSaggregation}
 \centering
  \begin{tabular}{r|c|ccccc}
   \hline\noalign{\smallskip}
   $\rho$ & $J_{n}$ & \multicolumn{4}{c}{$\delta_{i}$} \\
   \noalign{\smallskip}\hline
    16 & 0.9684 & 0.149 & 0.310 & 0.353 & 0.388 \\
    20 & 0.8941 & 0.057 & 0.146 & 0.243 & 0.308 \\
    50 & 0.8432 & 0.047 & 0.076 & 0.090 & 0.122 \\
   100 & 0.8376 & 0.018 & 0.030 & 0.037 & 0.041 \\
   16:128 & 0.8703 & $0.006$ & $0.006$ & $0.014$ & $0.014$ \\
   24:192 & 0.8611 & $0.005$ & $0.006$ & $0.006$ & $0.007$ \\
   \noalign{\smallskip}\hline
  \end{tabular}
\end{table}

\subsection{Results obtained with the KS function}
 \label{sSec:Aggregation functions}

We focus on the case $\overline{P_{c}} = 1.8$ and we solve again this by using the KS function to aggregate the eigenvalue constraints as in \eqref{eq:SingleConstraint}.

We first observe that, due to the conservative character of \eqref{eq:ksAggregation}, replacing \eqref{eq:EigConstraintMU} with \eqref{eq:SingleConstraint} amounts to solving a problem with a slightly higher buckling constraint. This can cause serious convergence issues, because of the delicate balance between stiffness and stability requirements. A first remedy to this is to scale the constraint \eqref{eq:SingleConstraint} according to the relative weight of $\mu_{1}$ in the aggregated measure. Therefore we define the scale factor
\begin{equation}
 \label{eq:ScaleFactor}
  s = \frac{\mu_{1}}{M\left[ \mu_{i} \right]}
\end{equation}
which is updated as the material penalization parameter $p$ is raised and jumps in the eigenvalues occur. This parameter is not considered in the sensitivity relation \eqref{eq:ksAggregationSens}, as it is updated only occasionally. This approach is conceptually similar to that adopted by \citet{le-etal_10a} al. for stress constrained problems.

The results obtained for some choices of the aggregation parameter $\rho$ are reported in \autoref{tab:ResultsKSaggregation} and the corresponding designs are shown in \autoref{fig:KSaggregationComparison}. From this figure we can also see the stiffness of the designs, compared to that of the design in \autoref{fig:OptimizedDesignsIncreasingPc} (g).

The aggregation parameter clearly plays an important role but, as a rule of thumb, designs reasonably close to those of \autoref{fig:OptimizedDesignsIncreasingPc} (g) are obtained when $\rho$ is set larger than 20. With $\rho = 50$ or above we obtain designs which also outperform the ones obtained with the original approach, even starting with a uniform material distribution. We recognize a two--fold effect of the aggregation parameter

\begin{enumerate}
 \item With a high value of $\rho$ we have an envelope which is closer to $\mu_{1}$, and therefore the overshooting of the constraint is reduced;
 \item The gaps between the eigenvalues will be reduced by increasing $\rho$ (see $\delta_{i}$ values in \autoref{tab:ResultsKSaggregation})
\end{enumerate}

The importance of choosing a relatively high value for the aggregation parameter, even for application to stress constraints, has been recently pointed out also by \citet{zhou-sigmund_17a}.

Regarding the second point, performing a continuation on the aggregation parameter gives very good results (see also \autoref{fig:KSaggregationComparison}). Here, the continuation strategy adopted was quite heuristic, as the parameter $\rho$ was simply doubled each 50 steps. For a more systematic continuation scheme, tuned on the progresses of the optimization, one can refer to \citet{poon-martins_07a}. However, we remark that also when employing continuation, the starting value of the aggregation parameter still has great importance. The choice of a low starting value (e.g. $\rho=16$) results in a bad design, whereas a very good one can be achieved by starting from $\rho = 24$ (see \autoref{fig:KSaggregationComparison}).

Similar considerations hold for the use of the $\rho$--norm function. From our experience a design equivalent to that of \autoref{fig:OptimizedDesignsIncreasingPc} (g) is obtained upon setting $\rho \geq 8$.

\section{Concluding discussion}
 \label{Sec:ConcludingDiscussion}

We have covered several issues related to topology optimization with linearized buckling constraints.

With the help of an illustrative 2D example we have demonstrated the influence of the buckling constraint on the optimized design. Furthermore, we have discussed the effect of some practices, as the use of non--conforming finite element approximations, the use of a simplified but inconsistent sensitivity expression, and the replacement of the separate eigenvalue constraints with an aggregated one.

The following main points are highlighted
\begin{itemize}
 \item We recognized the balance between stiffness and stability inherent to this problem, as well as the coalescence of several buckling modes contributing to the optimization. The implications of these features in complicating the optimization for high values of the buckling constraints have been discussed;
 \item Incompatible (or equivalently mixed) finite elements improve the accuracy of the buckling analysis on coarse meshes, compared to conforming ones. To the ends of the optimization, the use of these elements is advisable when the design is likely to consist of many thin features, which are prone to undergo local buckling;
 \item The use of inconsistent but simpler sensitivities has been shown to potentially lead to completely wrong results, and therefore should be avoided;
 \item Aggregation functions have proven to be very effective, leading to results which are essentially equivalent to those obtained by directly constraining the minimum eigenvalue. However, the aggregation parameter should be selected sufficiently high to obtain a close approximation to the extremal value, and not to overshoot too much the constraint;
\end{itemize}

As a final remark we conclude that several issues concerning TO with stability are far from being resolved, and a proper treatment of buckling, even in its linearized form, is relatively complicated. An issue not discussed here is the large computational cost due to the eigenvalue analyses, and this will be a topic for future investigations, aimed at the application and extension of efficient solution methods \citet{dunning-etal_16a,ferrari-etal_18a}.

The present paper should be a helpful discussion for such further researches, aimed at including buckling in large scale TO.

\begin{acknowledgements}
 The current project is supported by the Villum Fonden through the Villum Investigator Project ``InnoTop''. The authors are grateful to Prof. Pauli Pedersen for several fruitful discussions on the topic of the paper.
\end{acknowledgements}

\begin{small}
\bibliographystyle{spbasic}      
\bibliography{Database.bib}

\begin{thebibliography}{72}
\providecommand{\natexlab}[1]{#1}
\providecommand{\url}[1]{{#1}}
\providecommand{\urlprefix}{URL }
\expandafter\ifx\csname urlstyle\endcsname\relax
  \providecommand{\doi}[1]{DOI~\discretionary{}{}{}#1}\else
  \providecommand{\doi}{DOI~\discretionary{}{}{}\begingroup
  \urlstyle{rm}\Url}\fi
\providecommand{\eprint}[2][]{\url{#2}}

\bibitem[{Aage et~al(2017)Aage, Andreassen, Lazarov, and
  Sigmund}]{aage-etal_17a}
Aage N, Andreassen E, Lazarov BS, Sigmund O (2017) Giga--voxel computational
  morphogenesis for structural design. Nature 550(7674):84--86

\bibitem[{Achtziger(1999)}]{achtziger_99a}
Achtziger W (1999) Local stability of trusses in the context of topology
  optimization, {P}art {I}: exact modelling. Structural Optimization
  17(235--246)

\bibitem[{Armand and Lodier(1978)}]{armand-lodier_78a}
Armand JL, Lodier B (1978) Optimal design of bending elements. International
  Journal for Numerical Methods in Engineering 13:373--384

\bibitem[{Bathe and Dvorkin(1983)}]{bathe-dvorkin_83a}
Bathe KJ, Dvorkin E (1983) On the automatic solution of nonlinear finite
  element equations. Computers \& Structures 17(5--6):871--879

\bibitem[{Bends\o{}e and Sigmund(2004)}]{book:bendsoe-sigmund_2004}
Bends\o{}e MP, Sigmund O (2004) Topology Optimization: Theory, Methods and
  Applications. Springer

\bibitem[{Berke(1970)}]{berke_70a}
Berke L (1970) An efficient approach to the minimum weight design of deflection
  limited structures. AFFDDL-TM-70-4

\bibitem[{Bian and Feng(2017)}]{bian-feng_17a}
Bian X, Feng Y (2017) Large--scale buckling--constrained topology optimization
  based on assembly--free finite element analysis. Advances in Mechanical
  Engineering 9(9):1--12

\bibitem[{Bochenek and Tajs-Zieli{\'{n}}ska(2015)}]{bochenek-tajs_15a}
Bochenek B, Tajs-Zieli{\'{n}}ska K (2015) Minimal compliance topologies for
  maximal buckling load of columns. Structural and Multidisciplinary
  Optimization 51(5):1149--1157

\bibitem[{Brantman(1977)}]{brantman_77a}
Brantman R (1977) On the use of linearized instability analyses to investigate
  the buckling of nonsymmetrical systems. Acta Mechanica 26(1):75--89

\bibitem[{Bruyneel et~al(2008)Bruyneel, Colson, and
  Remouchamps}]{bruyneel-etal_08a}
Bruyneel M, Colson B, Remouchamps A (2008) Discussion on some convergence
  problems in buckling optimisation. Structural and Multidisciplinary
  Optimization 35(2):181--186

\bibitem[{Chen et~al(2004)Chen, Qi, Qi, and Teo}]{chen-etal_04a}
Chen X, Qi H, Qi L, Teo KL (2004) Smooth convex approximation to the maximum
  eigenvalue function. Journal of Global Optimization 30(2):253--270

\bibitem[{Cheng and Xu(2016)}]{cheng-xu_16}
Cheng G, Xu L (2016) Two--scale topology design optimization of stiffened or
  porous plate subject to out--of--plane buckling constraint. Structural and
  Multidisciplinary Optimization 54(5):1283--1296

\bibitem[{Chin and Kennedy(2016)}]{chin-kennedy_16a}
Chin TW, Kennedy GJ (2016) Large--scale compliance--minimization and buckling
  topology optimization of the undeformed common research model wing. In: AIAA
  SciTechForum

\bibitem[{Cook et~al(2001)Cook, Malkus, Plesha, and Witt}]{book:cook-etal01}
Cook RD, Malkus DS, Plesha ME, Witt RJ (2001) Concepts and Applications of
  Finite Element Analysis, 4th edn. Wiley

\bibitem[{Cox and McCarthy(1998)}]{cox-mccarthy_98a}
Cox S, McCarthy C (1998) The shape of the tallest column. SIAM J Math Anal
  29:547--554

\bibitem[{Cox and Overton(1992)}]{cox-overton_92a}
Cox S, Overton M (1992) On the optimal design of columns against buckling. SIAM
  J Math Anal 23:287--325

\bibitem[{Dunning et~al(2016)Dunning, Ovtchinnikov, Scott, and
  Kim}]{dunning-etal_16a}
Dunning PD, Ovtchinnikov E, Scott J, Kim A (2016) Level--set topology
  optimization with many linear buckling constraints using and efficient and
  robust eigensolver. International Journal for Numerical Methods in
  Engineering

\bibitem[{Duysinx and Sigmund(1998)}]{duysinx-sigmund_98a}
Duysinx P, Sigmund O (1998) New Developments in Handling Optimal Stress
  Constraints in Optimal Material Distributions, pp 1501--1509

\bibitem[{Ferrari et~al(2018)Ferrari, Lazarov, and Sigmund}]{ferrari-etal_18a}
Ferrari F, Lazarov BS, Sigmund O (2018) Eigenvalue topology optimization via
  efficient multilevel solution of the {F}requency {R}esponse. International
  Journal for Numerical Methods in Engineering 115(7):872--892

\bibitem[{Folgado and Rodrigues(1998)}]{folgado-rodrigues_98a}
Folgado J, Rodrigues H (1998) Structural optimization with a non-smooth
  buckling load criterion. Control \& Cybernetics 27(235--253)

\bibitem[{Frauenthal(1972)}]{frauenthal_72a}
Frauenthal JC (1972) Constrained optimal design of circular plates against
  buckling. Journal of Structural Mechanics 1:159--186

\bibitem[{Fr\"{o}ier et~al(1974)Fr\"{o}ier, Nilsson, and
  Samuelsson}]{froier-etal_74a}
Fr\"{o}ier M, Nilsson L, Samuelsson A (1974) The rectangular plane stress
  element by {T}urner, {P}ian and {W}ilson. International Journal for Numerical
  Methods in Engineering 8:433--437

\bibitem[{Gao and Ma(2015)}]{gao-ma_15a}
Gao X, Ma H (2015) Topology optimization of continuum structures under buckling
  constraints. Computers \& Structures 157:142--152

\bibitem[{Gao et~al(2017)Gao, Li, and Ma}]{gao-etal_17a}
Gao X, Li L, Ma H (2017) An adaptive continuation method for topology
  optimization of continuum structures considering buckling constraints.
  International Journal of applied mathematics 9(7):24

\bibitem[{Gravesen et~al(2011)Gravesen, Evgrafov, and
  Nguyen}]{gravesen-etal_11a}
Gravesen J, Evgrafov A, Nguyen DM (2011) On the sensitivities of multiple
  eigenvalues. Structural and Multidisciplinary Optimization 44(4):583--587

\bibitem[{Haftka and Gurdal(2012)}]{book:haftka2012}
Haftka R, Gurdal Z (2012) Elements of Structural Optimization. Solid Mechanics
  and Its Applications, Springer Netherlands

\bibitem[{Haftka and Prasad(1981)}]{haftka-prasad_81a}
Haftka R, Prasad B (1981) Optimum structural with plate bending elements -- {A}
  survey. AIAA Journal 19:517--522

\bibitem[{Hall et~al(1988)Hall, Cameron, and Grierson}]{hall-etal_88a}
Hall SK, Cameron GE, Grierson DE (1988) Least--weight design of steel
  frameframe accounting for {$P-\Delta$} effects. Structural Engineering {ASCE}
  115(6):1463--1475

\bibitem[{Jog and Haber(1996)}]{jog-haber_96a}
Jog CS, Haber RB (1996) Stability of finite element models for
  distributed-parameter optimization and topology design. Computer Methods in
  Applied Mechanics and Engineering 130(3):203 -- 226

\bibitem[{Kemmler et~al(2005)Kemmler, Lipka, and Ramm}]{kemmler-etal_05a}
Kemmler R, Lipka A, Ramm E (2005) Large deformations and stability in topology
  optimization. Structural and Multidisciplinary Optimization 30:459--476

\bibitem[{Kennedy and Hicken(2015)}]{kennedy-hicken_15a}
Kennedy GJ, Hicken JE (2015) Improved constraint--aggregation methods. Computer
  Methods in Applied Mechanics and Engineering 289(Supplement C):332 -- 354

\bibitem[{Kerr and T.~Soifer(1968)}]{kerr-soifer_68a}
Kerr AD, T~Soifer MT (1968) The linearization of the prebuckling state and its
  effect on the determined instability loads. Journal of Applied Mechanics 36

\bibitem[{Khot et~al(1976)Khot, Venkayya, and Berke}]{khot-etal_76a}
Khot NS, Venkayya VB, Berke L (1976) Optimum structural design with stability
  constraints. International Journal for Numerical Methods in Engineering
  10(5):1097--1114

\bibitem[{Kirmser and Hu(1995)}]{kirmser-hu_95a}
Kirmser PG, Hu KK (1995) The shape of the ideal column reconsidered. The
  Mathematical Intelligencer 15(3):62--67

\bibitem[{Kreisselmeier and Steinhauser(1979)}]{kreisselmeier-steinhauser_79a}
Kreisselmeier G, Steinhauser R (1979) Systematic control design by optimizing a
  vector performance index. IFAC Proceedings Volumes 12(7):113 -- 117, iFAC
  Symposium on computer Aided Design of Control Systems, Zurich, Switzerland,
  29-31 August

\bibitem[{Le et~al(2010)Le, Norato, Bruns, Ha, and Tortorelli}]{le-etal_10a}
Le C, Norato J, Bruns T, Ha C, Tortorelli D (2010) Stress-based topology
  optimization for continua. Structural and Multidisciplinary Optimization
  41(4):605--620

\bibitem[{Lee et~al(2012)Lee, James, and Martins}]{lee-etal_12a}
Lee E, James K, Martins J (2012) Stress--constrained topology optimization with
  design dependent loads. Structural and Multidisciplinary Optimization

\bibitem[{Lee et~al(2016)Lee, Ahn, and Yoo}]{lee-etal_16a}
Lee K, Ahn K, Yoo J (2016) A novel $p$--norm correction method for lightweight
  topology optimization under maximum stress constraints. Computers \&
  Structures 171(Supplement C):18 -- 30

\bibitem[{Lindgaard and Dahl(2013)}]{lindgaard-dahl_13a}
Lindgaard E, Dahl J (2013) On compliance and buckling objective functions in
  topology optimization of snap--through problems. Structural and
  Multidisciplinary Optimization 47:409--421

\bibitem[{Lund(2009)}]{lund_09a}
Lund E (2009) Buckling topology optimization of laminated multi--material
  composite shell structures. Composite Structures 91(2):158--167

\bibitem[{Manh et~al(2011)Manh, Evgrafov, Gersborg, and
  Gravesen}]{manh-etal_11a}
Manh ND, Evgrafov A, Gersborg AR, Gravesen J (2011) Isogeometric shape
  optimization of vibrating membranes. Computer Methods in Applied Mechanics
  and Engineering 200(13):1343--1353

\bibitem[{Mateus et~al(1997)Mateus, Mota~Soares, and
  Mota~Soares}]{mateus-etal_97a}
Mateus HC, Mota~Soares CM, Mota~Soares CA (1997) Buckling sensitivity analysis
  and optimal desing of thin laminated structure. Computer \& Structures
  64(1--4):461--472

\bibitem[{Munk et~al(2017)Munk, Vio, and Steven}]{munk-etal_16a}
Munk DJ, Vio GA, Steven GP (2017) A simple alternative formulation for
  structural optimisation with dynamic and buckling objectives. Structural and
  Multidisciplinary Optimization 55(3):969--986

\bibitem[{Neves et~al(1995)Neves, Rodrigues, and Guedes}]{neves-etal_95a}
Neves MM, Rodrigues H, Guedes JM (1995) Generalized topology design of
  structures with a buckling load criterion. Structural optimization
  10(2):71--78

\bibitem[{Neves et~al(2002)Neves, Sigmund, and Bends\o{}e}]{neves-etal_02a}
Neves MM, Sigmund O, Bends\o{}e MP (2002) Topology optimization of periodic
  microstructures with a penalization of highly localized buckling modes.
  International Journal for Numerical Methods in Engineering 54(6):809--834

\bibitem[{Ohsaki and Ikeda(2007)}]{book:makoto-ohsaki2007}
Ohsaki M, Ikeda K (2007) Stability and optimization of structures:
  {G}eneralized sensitivity analysis. Mechanical Engineering Series, Springer

\bibitem[{Olhoff and Rasmussen(1977)}]{olhoff-rasmussen_77a}
Olhoff N, Rasmussen SH (1977) On single and bimodal optimum buckling loads of
  clamped columns. International Journal of Solids and Structures
  13(7):605--614

\bibitem[{Pedersen and Pedersen(2018)}]{pedersen-pedersen_18b}
Pedersen NL, Pedersen P (2018) Buckling load optimization for 2{D} continuum
  models with alternative formulation for buckling load estimation. Structural
  and Multidisciplinary Optimization

\bibitem[{Pian(1964)}]{pian_64a}
Pian THH (1964) Derivation of element stiffness matrices by assumed stress
  distributions. AIAA Journal 2:1333--1336

\bibitem[{Pian and Sumihara(1984)}]{pian-sumihara_84a}
Pian THH, Sumihara K (1984) Rational approach for assumed stress finite
  elements. International Journal for Numerical Methods in Engineering
  20(9):1685--1695

\bibitem[{Poon and Martins(2007)}]{poon-martins_07a}
Poon NMK, Martins JRA (2007) An adaptive approach to contstraint aggregation
  using adjoint sensitivity analysis. Structural and Multidisciplinary
  Optimization 34:61--73

\bibitem[{Rahmatalla and Swan(2003)}]{rahmatalla-swan_03a}
Rahmatalla S, Swan C (2003) Continuum topology optimization of
  buckling--sensitive structures. AIAA Journal 41(6):1180--1189

\bibitem[{Raspanti et~al(2000)Raspanti, Bandoni, and
  Biegler}]{raspanti-etal_00a}
Raspanti CG, Bandoni JA, Biegler LT (2000) New strategies for flexibility
  analysis and desing under uncertainties. Computers \& Chemical Engineering
  24:2193--2209

\bibitem[{Rodrigues et~al(1995)Rodrigues, Guedes, and
  Bends\o{}e}]{rodrigues-etal_95a}
Rodrigues HC, Guedes JM, Bends\o{}e MP (1995) Necessary conditions for optimal
  design of structures with a nonsmooth eigenvalue based criterion. Structural
  Optimization 9:52--56

\bibitem[{Rojas-Labanda and Stolpe(2015)}]{labanda-stolpe_15a}
Rojas-Labanda S, Stolpe M (2015) Automatic penalty continuation in structural
  topology optimization. Structural and Multidisciplinary Optimization
  52:1205--1221

\bibitem[{Rozvany(1996)}]{rozvany_96a}
Rozvany G (1996) Difficulties in topology optimization with stress, local
  buckling and system stability constraints. Structural Optimization
  11:213--217

\bibitem[{Seyranian et~al(1994)Seyranian, Lund, and Olhoff}]{seyranian_94a}
Seyranian AP, Lund E, Olhoff N (1994) Multiple eigenvalues in structural
  optimization problems. Structural optimization 8(4):207--227

\bibitem[{Sigmund(2007)}]{sigmund_07a}
Sigmund O (2007) Morphology--based black and white filters for topology
  optimization. Structural and Multidisciplinary Optimization 33(4):401--424

\bibitem[{Simitses(1973)}]{simitses_73a}
Simitses GJ (1973) Optimal versus the stiffened circular plate. AIAA Journal
  11:1409--1412

\bibitem[{Simo and Rifai(1990)}]{simo-rifai_90a}
Simo JC, Rifai MS (1990) A class of mixed assumed strain methods and the method
  of incompatible modes. International Journal for Numerical Methods in
  Engineering 29:1595--1638

\bibitem[{Svanberg(1987)}]{svanberg_87a}
Svanberg K (1987) The method of moving asymptotes - {A} new method for
  structural optimization. International Journal for Numerical Methods in
  Engineering 24(2):359--373

\bibitem[{Szyszkowski and Watson(1988)}]{szyszkowsky-watson_88a}
Szyszkowski W, Watson L (1988) Optimization of the buckling load of columns and
  frames. Engineering Structures 10(4):249--256

\bibitem[{Thomsen et~al(2018)Thomsen, Wang, and Sigmund}]{thomsen-etal_18a}
Thomsen CR, Wang F, Sigmund O (2018) Buckling strength topology optimization of
  2{D} periodic materials based on linearized bifurcation analysis. Computer
  Methods in Applied Mechanics and Engineering 339:115--136

\bibitem[{Torii and Faria(2017)}]{torii-faria_17a}
Torii AJ, Faria JR (2017) Structural optimization considering smallest
  magnitude eigenvalues: a smooth approximation. Journal of the Brazilian
  Society of Mechanical Sciences and Engineering 39(5):1745--1754

\bibitem[{Turner et~al(1956)Turner, Clough, Martin, and Topp}]{turner-etal_56a}
Turner MJ, Clough RJ, Martin HC, Topp LJ (1956) Stiffness and deflection
  analysis of complex structures. Journal of Aerosol Science 23:805--823

\bibitem[{Verbart et~al(2017)Verbart, Langelaar, and van
  Keulen}]{verbart-etal_17a}
Verbart A, Langelaar M, van Keulen F (2017) A unified aggregation and
  relaxation approach for stress--constrained topology optimization. Structural
  and Multidisciplinary Optimization 55:663--679

\bibitem[{Wang et~al(2011)Wang, Lazarov, and Sigmund}]{wang-etal_11a}
Wang F, Lazarov B, Sigmund O (2011) On projection methods, convergence and
  robust formulations in topology optimization. Structural and
  Multidisciplinary Optimization 43(6):767--784

\bibitem[{Wilson et~al(1973)Wilson, Taylor, Doherty, and
  Glaboussi}]{wilson-etal_73a}
Wilson EL, Taylor RL, Doherty W, Glaboussi J (1973) Incompatible displacement
  models, Academic Press, pp 41--57

\bibitem[{Wu and Arora(1988)}]{wu-arora_1988}
Wu CC, Arora JS (1988) Design sensitivity analysis of non--linear buckling
  load. Computational Mechanics 3:129--140

\bibitem[{Yang and Chen(1996)}]{yang-chen_96a}
Yang RJ, Chen CJ (1996) Stress--based topology optimization. Structural
  Optimization 12(2):98--105

\bibitem[{Zhou(2004)}]{zhou_04a}
Zhou M (2004) Topology optimization of shell structures with linear buckling
  responses. In: {WCCM} {VI} in {B}eijing, {C}hina, pp 795--800

\bibitem[{Zhou and Sigmund(2017)}]{zhou-sigmund_17a}
Zhou M, Sigmund O (2017) On fully stressed design and $p$--norm measures in
  structural optimization. Structural and Multidisciplinary Optimization
  56(731--736)

\end{thebibliography}
\end{small}
\end{document}